\documentclass[11pt]{article}  
\usepackage{amsmath}
\usepackage{amssymb}
\usepackage{theorem}
\usepackage{euscript}
\usepackage{pstricks}
\newtheorem{theorem}{Theorem}[section]
\newtheorem{lemma}{Lemma}[section]

\numberwithin{equation}{section}
\title{\sffamily B-spline quasi-interpolant  representations  
and sampling recovery of functions with mixed smoothness
\author{ 
Dinh D\~ung \\[5mm]
Vietnam National University, Hanoi, Information Technology Institute \\
144 Xuan Thuy, Cau Giay, Hanoi, Vietnam\\
{\ttfamily dinhdung@vnu.edu.vn}\\[4mm]}}
\date{\ttfamily August 17, 2010 -- Version 0.95}
\def\II{{\mathbb I}}

\def\ZZ{{\mathbb Z}}
\def\NN{{\mathbb N}}
\def\RR{{\mathbb R}}
\def\Pp{{\mathcal P}}

\def\Ff{{\mathcal F}}
\def\MB{MB^\alpha_{p,\theta}}

\begin{document}
\maketitle

\begin{abstract}
Let $\xi = \{x^j\}_{j=1}^n$ be a grid of $n$ points 
in the $d$-cube ${\II}^d:=[0,1]^d$, and $\Phi = \{\varphi_j\}_{j =1}^n$ a family of  
$n$ functions on ${\II}^d$. 
We define the linear sampling algorithm $L_n(\Phi,\xi,\cdot)$ 
for an approximate recovery of a continuous function $f$ on ${\II}^d$ from the sampled values 
$f(x^1), ..., f(x^n)$, by 
\begin{equation*} 
L_n(\Phi,\xi,f) 
\ := \ 
\sum_{j=1}^n f(x^j)\varphi_j.
\end{equation*}
 For the Besov class $B^\alpha_{p,\theta}$ of mixed smoothness $\alpha$
(defined as the unit ball of the Besov space $\MB$), 
to study  optimality of $L_n(\Phi,\xi,\cdot)$ in $L_q({\II}^d)$
we use the quantity
\begin{equation*} 
r_n(B^\alpha_{p,\theta})_q 
\ := \ \inf_{H,\xi} \  \sup_{f \in B^\alpha_{p,\theta}} \, \|f - L_n(\Phi,\xi,f)\|_q, 
\end{equation*}
where the infimum is taken over all grids $\xi = \{x^j\}_{j=1}^n$ 
and all families $\Phi = \{\varphi_j\}_{j=1}^n$ in $L_q({\II}^d)$.
We  explicitly constructed linear sampling algorithms 
$L_n(\Phi,\xi,\cdot)$  on the grid 
$\xi = \ G^d(m):= \{(2^{-k_1}s_1,...,2^{-k_d}s_d) \in \II^d : \ k_1 + ... + k_d \le m\}$, 
with $\Phi$
a family of linear combinations of mixed B-splines which are mixed tensor products of either integer or 
half integer translated dilations of the centered B-spline of order $r$. 
The grid $G^d(m)$ is of the size $2^m m^{d-1}$ and sparse in comparing 
with the generating dyadic coordinate cube grid of the size $2^{dm}$. 
For various $0<p,q,\theta \le \infty$ and $1/p < \alpha < r$,
we proved upper bounds for the worst case error
$ \sup_{f \in B^\alpha_{p,\theta}} \, \|f - L_n(\Phi,\xi,f)\|_q$ which coincide with 
the asymptotic order of $r_n(B^\alpha_{p,\theta})_q$ in some cases.  
A key role in constructing these linear sampling algorithms, plays 
a quasi-interpolant representation of functions $f \in B^\alpha_{p,\theta}$ by mixed B-spline series
with the coefficient functionals which are explicitly constructed as
linear combinations of an absolute constant number of values of functions. Moreover, we proved 
that the quasi-norm of the Besov space $\MB$ is equivalent to a discrete quasi-norm in terms of  
the coefficient functionals.


  
\end{abstract}

\section{Introduction} 

The aim of the present paper is to investigate  linear  
sampling algorithms for recovery of functions on the unit $d$-cube 
${\II}^d:= [0,1]^d$, having a mixed smoothness. 
Let $\xi = \{x^j\}_{j=1}^n$ be a grid of $n$ points 
in ${\II}^d$, and $\Phi = \{\varphi_j\}_{j =1}^n$ a family of  
$n$ functions on ${\II}^d$. 
Then for  a continuous function $f$ on ${\II}^d$, we can define the linear sampling algorithm 
$L_n = L_n(\Phi,\xi,\cdot)$ 
for approximate recovering $f$ from the sampled values $f(x^1),..., f(x^n)$, by 
\begin{equation} \label{def:L_n(f)}
L_n(f) \ = \ L_n(\Phi,\xi,f) 
:= \ \sum_{j=1}^n f(x^j)\varphi_j.
\end{equation}
Let $L_q :=L_q({\II}^d), \ 0 < q \le \infty,$
denote the quasi-normed space 
of functions on ${\II}^d$ with the $q$th integral quasi-norm 
$\|\cdot\|_q$ for $0 < q < \infty,$ and 
the ess sup-norm $\|\cdot\|_{\infty}$ for $q = \infty$. 
The recovery error will be measured by $\|f - L_n(\Phi,\xi,f)\|_q$.

If $W$ is a class of continuous functions, 
$ \sup_{f \in B^\alpha_{p,\theta}} \, \|f - L_n(\Phi,\xi,f)\|_q$ is
 the worst case error of the recovery of functions $f$ from
$W$ by the linear sampling algorithm $L_n(\Phi,\xi,\cdot)$.
To study  optimality of linear sampling algorithms 
of the form \eqref{def:L_n(f)} for 
recovering $f \in W$ from $n$ their values,  we will use the quantity
\begin{equation} \label{def:r_n}
r_n(W)_q 
\ := \ \inf_{\xi, \Phi} \  \sup_{f \in W} \, \|f - L_n(\Phi,\xi,f)\|_q, 
\end{equation}
where the infimum is taken over all grids $\xi = \{x^j\}_{j=1}^n$ 
and all families $\Phi = \{\varphi_j\}_{j=1}^n$ in $L_q$.

A challenging problem in linear
sampling recovery of functions from a class $W$ with a given mixed smoothness, is to construct 
a sampling algorithm $L_n(\Phi,\xi,\cdot)$ with 
an appropriate  sampling grid  $\xi = \{x^j\}_{j=1}^n$ and family $\Phi = \{\varphi_j\}_{j=1}^n$
which would be asymptotically optimal in terms of the quantity $r_n(W)_q$. 

For periodic functions Smolyak \cite{S} first constructed 
a specific linear sampling algorithm based on the de la Vallee Poussin kernel 
and the following dyadic grid in $\II^d$
\begin{equation*} 
G^d(m)
:= \ 
\{(2^{-k_1}s_1,...,2^{-k_d}s_d) \in \II^d : \  k \in \Delta(m)\}
\ = \
\{ 2^{-k}s: k \in \Delta(m),\ s \in I^d(k)\}.
\end{equation*}
Here and in what follows, we use the notations:
$xy := (x_1y_1,..., x_dy_d)$; 
$2^x := (2^{x_1},...,2^{x_d})$;
$|x|_1 := \sum_{i=1}^d |x_i|$ for $x, y \in {\RR}^d$;
$\Delta(m) := \{k \in {\ZZ}^d_+: |k|_1 \le m\}$;  
$I^d(k):= \{s \in {\ZZ}^d_+: 0 \le s_i \le 2^{k_i}, \ i \in N[d]\}$; 
$N[d]$ denotes the set of all natural numbers from $1$ to $d$; $x_i$ denotes the $i$th coordinate 
of $x \in \RR^d$, i.e., $x := (x_1,..., x_d)$.
Temlyakov \cite{Te1}, \cite{Te2}, \cite{Te3} and Dinh Dung \cite{Di2}--\cite{Di4} developed
Smolyak's construction for study the asymptotic order of $r_n(W)_q$ for  
periodic Sobolev classes $W^\alpha_p$ and H\"older classes $H^\alpha_p$ as well their intersection.
In particular, the first asymptotic order 
\begin{equation*} 
r_n(H^\alpha_p)_q
 \ \asymp \ 
(n^{-1} \log^{d-1}n)^{\alpha - 1/p + 1/q}(\log^{d-1}n)^{1/q }, \ 1< p < q \le 2, \ \alpha > 1/p,
\end{equation*}
was obtained in \cite{Di2}--\cite{Di3}. 
For non-periodic functions of mixed smoothness $1/p < \alpha \le 2$,
this problem has been recently studied by Sickel and Ullrich \cite{SU},
using the mixed tensor product of piecewise linear B-splines (of order $2$) and the grid $G^d(m)$.
It is interesting to notice that the linear sampling algorithms 
considered by above mentioned authors
are interpolating at the grid $G^d(m)$.

Naturally, the quantity $r_n(W)_q $ of optimal linear sampling recovery is related 
to the problem of optimal linear approximation in terms of
the linear $n$-width $\lambda_n(W)_q $ introduced by Tikhomirov \cite{Ti}:
\begin{equation*} 
\lambda_n(W)_q 
\ := \ \inf_{A_n} \  \sup_{f \in W} \, \|f - A_n(f)\|_q, 
\end{equation*}
where the infimum is taken over all linear operators $A_n$ of rank $n$ in $L_q$.
The linear $n$-width $\lambda_n(W)_q $
was studied in \cite{G}, \cite{R1}, \cite{R2}, ect. for 
 various classes $W$ of functions with mixed smoothness.  The inequality $r_n \ge \lambda_n$ is 
quite useful in investigation of the (asymptotic) optimality of a given
linear sampling algorithm. It also allows to establish a lower bound of $r_n$ via 
a known lower bound of $\lambda_n$. 

In the present paper, we continue to research this problem. We will
 take functions to be recovered from the Besov class
$B^\alpha_{p,\theta}$ of functions on ${\II}^d$, which is defined as the unit ball of the Besov space 
$\MB$  having mixed smoothness $\alpha$. 
For functions in $B^\alpha_{p,\theta}$, we will construct linear sampling algorithms 
$L_n(\Phi,\xi,\cdot)$ on the grid $\xi=G^d(m)$ with $\Phi$
a family of linear combinations of mixed B-splines which are mixed tensor products of either integer or 
half integer translated dilations of the centered B-spline of order $r > \alpha$. We will 
be concerned with the worst case error of the recovery of $B^\alpha_{p,\theta}$ in the space $L_q$ by 
these linear sampling algorithms and their asymptotic optimality in terms of the quantity 
$r_n(B^\alpha_{p,\theta})_q$ for various $0 < p, q, \theta \le \infty$ and $1/p \le \alpha < r$. 
A key role in constructing these linear sampling algorithms, plays 
a quasi-interpolant representation of functions $f \in \MB$ by mixed 
B-spline series which will be explicitly constructed. Let us give a sketch of the 
main results of the present paper. 

We first describe representations by mixed B-spline series constructed on the basic of quasi-interpolants.
 For a given natural number $r,$ let  $M$ be the   
centered B-spline of order $r$ with support $[-r/2,r/2]$ and 
knots at the points $-r/2,-r/2 + 1,...,r/2 - 1, r/2,$. 
We define the integer translated dilation $M_{k,s}$ of $M$ by   
\begin{equation*}
M_{k,s}(x):= \ M(2^k x - s), \ k \in {\ZZ}_+, \ s \in \ZZ, 
\end{equation*}
and the mixed $d$-variable B-spline $M_{k,s}$ by
\begin{equation} \label{def:Mixed[M_{k,s}]}
M_{k,s}(x):=  \ \prod_{i=1}^d M_{k_i, s_i}( x_i),  
\ k \in {\ZZ}^d_+, \ s \in {\ZZ}^d,
\end{equation}
where ${\ZZ}_+$ is the set of all non-negative integers,  
${\ZZ}^d_+:= \{s \in {\ZZ}^d: s_i \ge 0, \ i \in N[d] \}$.
Further, we define 
the half integer translated dilation $M^*_{k,s}$ of $M$ by     
\begin{equation*}
M^*_{k,s}(x):= \ M(2^k x - s/2), \ k \in {\ZZ}_+, \ s \in \ZZ, 
\end{equation*}
and the mixed $d$-variable B-spline $M^*_{k,s}$ by
\begin{equation*} 
M^*_{k,s}(x):= \ \ \prod_{i=1}^d M^*_{k_i, s_i}( x_i), 
\ k \in {\ZZ}^d_+, \ s \in {\ZZ}^d.
\end{equation*}
In what follows, the B-spline $M$ will be fixed. We will denote
$M^r_{k,s}:= M_{k,s}$ if the order $r$ of $M$ is even, 
and $M^r_{k,s}:= M^*_{k,s}$ if the order $r$ of $M$  is odd.

Let $0 < p, \theta \le \infty,$ and $1/p <\alpha < \min (r, r - 1 + 1/p)$. 
Then we prove the following mixed B-spline quasi-interpolant representation of 
 functions $f \in \MB$. Namely, a function $f$ in the Besov space  
$\MB$ can be represented by the mixed B-spline series 
\begin{equation} \label{Representation1}
f \ =  \
\sum_{k \in {\ZZ}^d_+} \sum_{s \in J_r^d(k)}c^r_{k,s}(f)M^r_{k,s},
\end{equation}
converging in the quasi-norm of $\MB$, where $J_r^d(k)$ is the set of $s$ for which $M^r_{k,s}$
do not vanish identically on  ${\II}^d$, and 
the coefficient functionals $c^r_{k,s}(f)$ explicitly constructed as
linear combinations of an absolute constant number of values of $f$ which does not depend 
on neither $k,s$ nor $f$. 
Moreover, we prove that the quasi-norm of $\MB$  
is equivalent to some discrete quasi-norm in terms of the coefficient functionals 
$c^r_{k,s}(f)$. B-spline quasi-interpolant representations of 
 functions from the isotropic Besov sapces has been constructed in \cite{Di6}, \cite{Di7}.
Different B-spline quasi-interpolant representations were considered in \cite{DP}.
Both these representations were constructed on the basic
of B-spline quasi-intepolants. The reader can see the books \cite{C}, \cite{BHR} 
for survey and details on quasi-interpolants. 

Let us construct linear sampling algorithms 
$L_n(\Phi,\xi,\cdot)$ on the grid $\xi=G^d(m)$ on the basic 
of the representation \eqref{Representation1}. 
For $m \in {\ZZ}_+$, let the linear operator $R_m$ be defined 
for functions $f$ on ${\II}^d$ by
\begin{equation} \label{def:R_m}
R_m(f) 
\ := \ 
\sum_{k \in \Delta(m)} \ \sum_{s \in J_r^d(k)}c^r_{k,s}(f)M^r_{k,s}. 
\end{equation} 
If ${\bar m}$ is the largest of $m$ such that 
\begin{equation*} 
2^m  m^{d-1} \asymp  |G^d(m)| \le n
\end{equation*} 
for a given $n$, where $|A|$ denotes the cardinality of $A$, then the operator $R_{{\bar m}}$
is a linear sampling algorithm 
of the form \eqref{def:L_n(f)} on the grid $G^d({\bar m})$. More precisely, 
\begin{equation*} 
R_{{\bar m}}(f) 
\ = \ 
L_n(\Phi,\xi,f) 
\ = \ 
\sum_{(k,s) \in G^d({\bar m})} f(2^{-k}s) \psi_{k,s}, 
\end{equation*} 
where $\psi_{k,s}$ are  explicitly constructed as linear combinations of an absolute constant 
 of  B-splines $M^r_{k,j}$, which does not depend 
on neither $k,s$ nor $f$. It is worth to emphasize  that 
the grid $G^d(m)$ is of the size $2^m m^{d-1}$ and sparse in comparing 
with the generating dyadic coordinate cube grid of the size $2^{dm}$. 
We give now a brief of our results concerning with the worst case error of 
the recovery of functions $f$ from $B^\alpha_{p,\theta}$ by the linear sampling algorithms
$R_{\bar m}(f)$ and their asymptotic optimality.

We use the notations: $x_+ := \max(0,x)$ for $x \in \RR;$  
$A_n(f) \ll B_n(f)$ if $A_n(f) \le CB_n(f)$ with 
$C$ an absolute constant not depending on $n$ and/or $f \in W,$ and 
$A_n(f) \asymp B_n(f)$ if $A_n(f) \ll B_n(f)$ and $B_n(f) \ll A_n(f).$ 
Let us introduce  the abbreviations: 
\begin{equation*} 
E(m):= \sup_{f \in B^\alpha_{p,\theta}}\|f - R_m(f)\|_q, \quad
 r_n := r_n(B^\alpha_{p,\theta})_q.
\end{equation*} 

Let $\ 0 < p, q, \theta \le \infty$ and $1/p < \alpha < r$. 
Then we have the following upper bound of $r_n$ and $E({\bar m})$. 
\begin{itemize}
\item[{\rm (i)}] For $p \ge q$,
\begin{equation} \label{Intr:[r_n<,p>q]}
r_n
\ \ll \
E({\bar m})
\ \ll \ 
\begin{cases}
(n^{-1} \log^{d-1}n)^\alpha, \ & \theta \le \min(q,1), \\
(n^{-1} \log^{d-1}n)^\alpha (\log^{d-1}n)^{ 1/q - 1/\theta}, \ & \theta > \min(q,1), \ q \le 1, \\
(n^{-1} \log^{d-1}n)^\alpha (\log^{d-1}n)^{ 1 - 1/\theta}, \ & \theta > \min(q,1), \ q > 1.
\end{cases}
\end{equation}
\item[{\rm (ii)}] For $p < q$, 
\begin{equation} \label{Intr:[r_n<,p<q]}
r_n
\ \ll \
E({\bar m})
 \ \ll \  
\begin{cases}
(n^{-1} \log^{d-1}n)^{\alpha - 1/p + 1/q}(\log^{d-1}n)^{(1/q - 1/\theta)_+}, \ & q < \infty, \\
(n^{-1} \log^{d-1}n)^{\alpha - 1/p}(\log^{d-1}n)^{(1 - 1/\theta)_+}, \ & q = \infty.
\end{cases}
\end{equation}
\end{itemize} 

From the embedding of $\MB$ into the isotropic Besov space of smootthness $d\alpha$ and known
asymptotic order of the quantity \eqref{def:r_n} of its unit ball in $L_q$
(see \cite{Di1}, \cite{K}, \cite{No}, \cite{NoT}, \cite{Te3}) it follows that for
$0 < p,q \le \infty, \ 0 < \theta \le \infty$ and $\alpha > 1/p$,
there always holds the lower bound 
$r_n \ \gg \ n^{- \alpha + (1/p - 1/q)_+}$. 
However, this estimation is too rough and does not lead 
to the asymptotic order. By use of the inequality $\lambda_n(B^\alpha_{p,\theta})_q \ge r_n$ and
known results on $\lambda_n(B^\alpha_{p,\theta})_q$ \cite{G}, \cite{R1},
from \eqref{Intr:[r_n<,p>q]} and \eqref{Intr:[r_n<,p<q]}
we obtain the asymptotic order of $r_n$ for some cases.
More precisely, we have the following asymptotic orders of $r_n$ and $E({\bar m})$ which show
the asymptotic optimality of the linear sampling algorithms $R_{{\bar m}}$. 
\begin{itemize}
\item[{\rm (i)}] For $p \ge q$ and $\theta \le 1$,
\begin{equation} \label{Intr:[r_n><,p>q]}
E({\bar m})
 \ \asymp \ 
r_n
 \ \asymp \ 
 (n^{-1} \log^{d-1}n)^\alpha, \ 
\begin{cases}
 2 \le q < p < \infty, \\
1 < p = q \le \infty.
\end{cases}
\end{equation}
\item[{\rm (ii)}] For $1 < p < q < \infty$, 
\begin{equation} \label{Intr:[r_n><,p<q]}
E({\bar m})
 \ \asymp \ 
r_n
 \ \asymp \ 
(n^{-1} \log^{d-1}n)^{\alpha - 1/p + 1/q}(\log^{d-1}n)^{(1/q - 1/\theta)_+}, \
 \begin{cases}
2 \le p, \ 2 \le \theta \le q,  \\
q \le 2.
\end{cases}
\end{equation}
\end{itemize}

The present paper is organized as follows. 
In Section \ref{Quasi-interpolant}, we give a necessary background of Besov spaces of mixed 
smoothness,  B-spline quasi-interpolants, and prove a theorem on the mixed B-spline 
quasi-iterpolant representation \eqref{Representation1}
and a relevant discrete equivalent quasi-norm  for the Besov space of mixed smoothness $\MB$. 
In Section \ref{Sampling recovery}, we prove the upper bounds 
\eqref{Intr:[r_n<,p>q]}--\eqref{Intr:[r_n<,p<q]} and the asymptotic orders 
\eqref{Intr:[r_n><,p>q]}--\eqref{Intr:[r_n><,p<q]}. In Section \ref{Faber-Schauder},
we consider interpolant representations by series of 
the mixed tensor product of piecewise constant or piecewise linear B-splines.
 In Section \ref{Appendix}, we present some auxiliary results.

\section{B-spline quasi-interpolant representations} 
\label{Quasi-interpolant}

Let us introduce 
Besov spaces of functions with mixed smoothness
and give necessary knowledge of them. 
For univariate functions the $l$th difference operator $\Delta_h^l$ is defined by 
\begin{equation*}
\Delta_h^lf(x) := \
\sum_{j =0}^l (-1)^{l - j} \binom{l}{j} f(x + jh).
\end{equation*}
If $e$ is any subset of $N[d]$, for multivariate functions 
the mixed $(l,e)$th difference operator $\Delta_h^{l,e}$ is defined by 
\begin{equation*}
\Delta_h^{l,e} := \
\prod_{i \in e} \Delta_{h_i}^l, \ \Delta_h^{l,\emptyset} = I,
\end{equation*}
where the univariate operator
$\Delta_{h_i}^l$ is applied to the univariate function $f$ by considering $f$ as a 
function of  variable $x_i$ with the other variables held fixed. 
For a domain $\Omega$ in ${\RR}^d$, denote by  
$L_p(\Omega)$ the quasi-normed space 
of functions on $\Omega$ with the $p$th integral quasi-norm 
$\|\cdot\|_{p,\Omega}$ for $0 < p < \infty,$ and 
the ess sup-norm $\|\cdot\|_{\infty,\Omega}$ for $p = \infty$.
Let
\begin{equation*}
\omega_l^e(f,t)_p:= \sup_{|h_i| < t_i, i \in e}\|\Delta_h^{l,e}(f)\|_{p,{\II}^d(h,e)}, \ t \in {\II}^d,
\end{equation*} 
be the mixed $(l,e)$th modulus of smoothness of $f$,
 where ${\II}^d(h,e):= \{ x \in {\II}^d : x_i, x_i + lh_i \in \II, \ i \in e \}$ (in particular,
 $\omega_l^{\emptyset}(f,t)_p = \|f\|_p$).
We will need the following modified $(l,e)$th mixed modulus of smoothness
\begin{equation*}
w_l^e(f,t)_p 
\ := \ 
\left( \prod_{i \in e} t_i^{- 1} \int_{U(t,e)} \int_{{\II}^d(h,e)}
|\Delta_h^l(f,x)|^p \ dx \ dh \right)^{1/p},
\end{equation*}
where $U(t,e):= \{ x \in {\II}^d: |x_i| \le t, \ i \in e\}$.
 There hold the following inequalities
\begin{equation} \label{ineq:w_l><omega_l}
C_1 w_l^e(f,t)_p \ \le \omega_l^e (f,t)_p \ \le \ C_2 w_l^e(f,t)_p 
\end{equation}
with constants $C_1,C_2$ which depend on $l,p,d$ only. A proof of these inequalities 
for the univariate modulus of smoothness is given in \cite{PP}.
They can be proven in a similar way for the multivariate $(l,e)$th mixed modulus of smoothness.  

If $0 <  p, \theta \le \infty$, 
$\alpha > 0$ and $l > \alpha$, 
we introduce the quasi-semi-norm 
$|f|_{B_{p, \theta}^{\alpha,e}}$ for functions $f \in L_p$ by
\begin{equation*} \label{BesovSeminorm}
|f|_{MB_{p, \theta}^{\alpha,e} }:= 
\begin{cases}
 \ \left(\int_{{\II}^d} \{ \prod_{i \in e} t_i^{-\alpha_i}
\omega_l^e(f,t)_p \}^ \theta \prod_{i \in e} t_i^{-1}dt \right)^{1/\theta}, 
& \theta < \infty, \\
  \sup_{t \in {\II}^d} \ \prod_{i \in e} t_i^{-\alpha_i}\omega_l^e(f,t)_p,  & \theta = \infty
\end{cases}
\end{equation*}
(in particular, $|f|_{MB_{p, \theta}^{\alpha,\emptyset}} = \|f\|_p$).

For $0 <  p, \theta \le \infty$ and $0 < \alpha < l,$ the Besov space 
$\MB$ is defined as the set of  functions $f \in L_p$ 
for which the Besov quasi-norm $\|f\|_{\MB}$ is finite. 
The  Besov quasi-norm is defined by
\begin{equation*} 
B(f) \ = \ \|f\|_{\MB}
:= \ 
 \sum_{e \subset N[d]} |f|_{MB_{p, \theta}^{\alpha,e} }.
\end{equation*}

We will study the linear sampling recovery of functions from the Besov class
\begin{equation*}
B^\alpha_{p,\theta}
:= \ 
\{f \in \MB: \ B(f)\le 1\},
\end{equation*}
 with the restriction on the smoothness $\alpha > 1/p$,
which provides the compact embedding of $\MB$ into $C({\II}^d)$, 
the space of continuous functions on $\II^d$ with max-norm.
 We will also study this problem for $B^\alpha_{p,\theta}$ with the restrictions
$\alpha = 1/p$ and $p \le \min(1, \theta)$ which is a sufficient condition for  
the continuous embedding of $\MB$ into $C({\II}^d)$.
In both these cases, $B^\alpha_{p,\theta}$ can be considered as a subset in $C(\II^d)$. 

For any $e \subset N[d]$, put 
${\ZZ}^d_+(e):= \{s \in {\ZZ}^d_+: s_i = 0 , \ i \notin e\}$
(in particular, ${\ZZ}^d_+(\emptyset)= \{0\}$ and ${\ZZ}^d_+(N[d])= {\ZZ}^d_+$).
If $\{g_k\}_{k \in {\ZZ}^d_+(e)}$ is a sequence 
whose component functions $g_k$ are in 
$L_p,$ for $0 < p, \theta \le \infty$ and $\beta \ge 0$ we define the 
$b_\theta^{\beta,e}(L_p)$ ``quasi-norms"
\begin{equation*}
\ \|\{ g_k\}\|_{b_\theta^{\beta,e}(L_p)}
\ := \ 
 \biggl(\sum_{k \in {\ZZ}^d_+(e)} 
 \left(2^{\beta |k|_1}\| g_k \|_p\right)^\theta \biggl)^{1/\theta} 
\end{equation*}
with the usual change to a supremum when $\theta = \infty.$ 
When $\{g_k\}_{k \in {\ZZ}^d_+(e)}$ is a positive sequence, we replace 
$\| g_k \|_p$ by $|g_k|$ and denote the corresponding quasi-norm by 
$\|\{g_k\}\|_{b_\theta^{\beta,e}}$. 

For the Besov space $\MB$, from the definition and properties of  the mixed $(l,e)$th 
modulus of smoothness it is easy to verify that 
there is the following quasi-norm equivalence
\begin{equation*} 
B(f) 
\ \asymp \ B_1(f)
\ := \
\sum_{e \subset N[d]} \|\{\omega_l^e(f,2^{-k})_p\}\|_{b_\theta^{\alpha,e}}. 
\end{equation*}

Let $\Lambda = \{\lambda(s)\}_{j \in P(\mu)}$ be a finite even sequence, i.e., 
$\lambda(-j) = \lambda(j),$ where $P(\mu):= \{j \in  \ZZ: \ |j| \le \mu \}$ 
and $\mu \ge r/2 - 1$. 
We define the linear operator $Q$ for functions $f$ on $\RR$ by  
\begin{equation} \label{def:Q}
Q(f,x):= \ \sum_{s \in \ZZ} \Lambda (f,s)M(x-s), 
\end{equation} 
where
\begin{equation} \label{def:Lambda}
\Lambda (f,s):= \ \sum_{j \in P(\mu)} \lambda (j) f(s-j).
\end{equation}
The operator $Q$ is bounded in $C(\RR)$ and 
\begin{equation*}
\|Q(f)\|_{C(\RR)} \le \|\Lambda \|\|f\|_{C(\RR)}  
\end{equation*}
for each $f \in C(\RR),$ where
\begin{equation*}
\|\Lambda \|= \ \sum_{j \in P(\mu)} |\lambda (j)|. 
\end{equation*}

Moreover, $Q$ is local in the following sense. There is 
a positive number $\delta > 0$ such that for  any  
$f \in C(\RR)$ and $x \in \RR,$ $Q(f,x)$ 
depends only on the value $f(y)$ at an absolute constant number of points 
$y$ with $|y - x|\le \delta$. 
We will require $Q$ to reproduce  the space 
$\Pp_{r-1}$ of polynomials of order at most $r - 1$ , that is,  
$Q(p) \ = \ p, \ p \in \Pp_{r-1}$. 
An operator $Q$ of the form \eqref{def:Q}--\eqref{def:Lambda} reproducing 
$\Pp_{r-1}$, is called a {\it quasi-interpolant in} $C(\RR).$ 

There are many ways to construct  quasi-interpolants.
 A method of construction via 
Neumann series was suggested by Chui and Diamond \cite{CD} 
(see also \cite[p. 100--109]{C}). A necessary and sufficient condition 
of reproducing $\Pp_{r-1}$ for operators $Q$ of the form \eqref{def:Q}--\eqref{def:Lambda} with even $r$ and 
$\mu \ge r/2$, was established 
in \cite{BSSV}. De Bore and Fix \cite{BF} 
introduced another quasi-interpolant 
based on the values of derivatives. 

Let us give some examples of quasi-interpolants. 
The simplest example is 
a piecewise constant quasi-interpolant which is defined for $r=1$ by
\begin{equation*} 
Q(f,x):= \ \sum_{s \in \ZZ} f(s) M(x-s), 
\end{equation*} 
where $M$ is the   
symmetric piecewise constant B-spline with support $[-1/2,1/2]$ and 
knots at the half integer points $-1/2, 1/2$. 
A piecewise linear quasi-interpolant is defined for $r=2$ by
\begin{equation} \label{def:Q,r=2}
Q(f,x):= \ \sum_{s \in \ZZ} f(s) M(x-s), 
\end{equation} 
where $M$ is the   
symmetric piecewise linear B-spline with support $[-1,1]$ and 
knots at the integer points $-1, 0, 1$.
This quasi-interpolant is also called {nodal} and directly related to the classical Faber-Schauder basic. 
We will revisit  it in Section \ref{Faber-Schauder}. A quadric quasi-interpolant is defined for $r=3$ by
\begin{equation*} 
Q(f,x):= \ \sum_{s \in \ZZ} \frac {1}{8} \{- f(s-1) + 10f(s) - f(s+1)\} M(x-s), 
\end{equation*} 
where $M$ is the symmetric quadric B-spline with support $[-3/2,3/2]$ and 
knots at the half integer points $-3/2, -1/2, 1/2, 3/2$.
Another example is the cubic quasi-interpolant defined for $r=4$ by
\begin{equation*} 
Q(f,x):= \ \sum_{s \in \ZZ} \frac {1}{6} \{- f(s-1) + 8f(s) - f(s+1)\} M(x-s), 
\end{equation*} 
where $M$ is the symmetric cubic B-spline with support $[-2,2]$ and 
knots at the integer points $-2, -1, 0, 1, 2$.

If $Q$ is a quasi-interpolant of the form 
\eqref{def:Q}--\eqref{def:Lambda}, for $h > 0$ and a function $f$ on $\RR$, 
we define the operator $Q(\cdot;h)$ by
\begin{equation*}
Q(f;h) = \ \sigma_h \circ Q \circ \sigma_{1/h}(f),
\end{equation*}
where $\sigma_h(f,x) = \ f(x/h)$.
By definition it is easy to see that 
\begin{equation*}
Q(f,x;h)= \ 
\sum_{k}\Lambda (f,k;h)M(h^{-1}x-k),
\end{equation*}
where 
\begin{equation*}
\Lambda (f,k;h):= 
 \ \sum_{j \in P(\mu)} \lambda (j) f(h(k-j)).
\end{equation*}

The operator $Q(\cdot;h)$ has the same properties as $Q$: it is a local bounded linear 
operator in $C(\RR)$ and reproduces the polynomials from $\Pp_{r-1}.$ 
Moreover, it gives a good approximation for smooth functions \cite[p. 63--65]{BHR}.
We will also call it a \emph{quasi-interpolant for} $C(\RR).$
However, the quasi-interpolant $Q(\cdot;h)$ is not defined for 
a function $f$ on $\II$, and therefore, not appropriate for an 
approximate sampling recovery of $f$ from its sampled values at points in $\II$. 
An approach to construct a quasi-interpolant for functions on $\II$ is to extend it 
by interpolation Lagrange polynomials. This approach has been proposed in 
\cite{Di6} for the univariate case. Let us recall it.

For a non-negative integer $k,$ 
we put $x_j = j2^{-k},  j \in \ZZ.$ If $f$ is a function on $\II,$ let  
\begin{equation*} 
\begin{aligned}
U_k(f,x) & := \ f(x_0) \ + \ \sum_{s=1}^{r-1}
\ \frac{2^{sk} \Delta_{2^{-k}}^sf(x_0)}{s!}
\ \prod_{j=0}^{s-1} (x - x_j), \\
V_k(f,x) & := \ f(x_{2^k - r + 1}) \ + \ \sum_{s=1}^{r-1}
\ \frac{2^{sk} \Delta_{2^{-k}}^sf(x_{2^k - r + 1})}{s!}
\ \prod_{j=0}^{s-1} (x - x_{2^k - r + 1 + j})
\end{aligned}
\end{equation*} 
be the $(r - 1)$th Lagrange polynomials interpolating $f$
at the $r$ left end points $x_0, x_1,..., x_{r-1},$ and 
$r$ right end points 
$x_{2^k - r + 1}, x_{2^k - r + 3},..., x_{2^k},$ 
of the interval $\II,$ respectively. 
The function $f_k$ is defined as an extension of 
$f$ on $\RR$ by the formula
\begin{equation*} 
f_k (x):= \
\begin{cases}
U_k(f,x), \ & x < 0, \\
f(x), \ & 0 \le x \le 1, \\
V_k(f,x), \ & x >1.
\end{cases}
\end{equation*} 
Obviously, if $f$ is continuous on $\II$, then $f_k$ is a continuous function on $\RR.$ 
Let $Q$ be a quasi-interpolant of the form \eqref{def:Q}--\eqref{def:Lambda} in $C({\RR}).$ 
Put ${\Bar {\ZZ}}_+ := \{k \in \ZZ: k \ge -1 \}$.
If $k \in {\Bar {\ZZ}}_+ $, we introduce the operator $Q_k$  by  
\begin{equation*}
Q_k(f,x) = \ Q(f_k,x;2^{-k}),\ \text{and} \ Q_{-1}(f,x):= 0, \  x \in \II, 
\end{equation*}
for a function $f$ on $\II$. We have for $k \in {\ZZ}_+$,
\begin{equation} \label{eq:Q_k}
Q_k(f,x)  \ = \ 
\sum_{s \in J(k)} a_{k,s}(f)M_{k,s}(x), \ \forall x \in \II, 
\end{equation}
where $J(k) := \ \{s \in \ZZ:\ -r/2 < s <   2^k + r/2 \}$ 
is the set of $s$ for which $M_{k,s}$ 
do not vanish identically on  $\II,$ and the coefficient functional $a_{k,s}$ is defined by
\begin{equation*} 
a_{k,s}(f):= \ \Lambda(f_k,s;{2^{-k}}) 
= \   
\sum_{|j| \le \mu} \lambda (j) f_k(2^{-k}(s-j)).
\end{equation*}

Put ${\Bar {\ZZ}}^d_+ := \{k \in {\ZZ}^d_+: k_i \ge -1, \ i \in N[d]\}$. 
For $k \in {\Bar {\ZZ}}^d_+$, let the mixed operator $Q_k$ be defined by
\begin{equation} \label{def:Mixed[Q_k]} 
Q_k:= \prod_{i=1}^d  Q_{k_i},
\end{equation}
where the univariate operator
$Q_{k_i}$ is applied to the univariate function $f$ by considering $f$ as a 
function of  variable $x_i$ with the other variables held fixed.

We have
\begin{equation*} 
Q_k(f,x)  \ = \ 
\sum_{s \in J^d(k)} a_{k,s}(f)M_{k,s}(x), \quad \forall x \in {\II}^d, 
\end{equation*}
where $M_{k,s}$ is the mixed B-spline  defined in \eqref{def:Mixed[M_{k,s}]},  
$J^d(k) := \ \{s \in {\ZZ}^d:\ \ - r/2 < s_i < 2^{k_i} + r/2, \ i \in N[d]\}$ 
is the set of $s$ for which $M_{k,s}$ do not vanish identically on  ${\II}^d$, 
\begin{equation} \label{def:Mixed[a_{k,s}(f)]}
a_{k,s}(f) 
\ := \   
a_{k_1,s_1}((a_{k_2,s_2}(...a_{k_d,s_d}(f))),
\end{equation}
and the univariate coefficient functional
$a_{k_i,s_i}$ is applied to the univariate function $f$ by considering $f$ as a 
function of  variable $x_i$ with the other variables held fixed. 

The operator $Q_k$ is a local bounded linear 
mapping in $C({\II}^d)$  and reproducing $\Pp_{r-1}^d$ the space 
of polynomials of order at most $r - 1$ in each variable $x_i$.   More precisely, there is 
a positive number $\delta > 0$ such that for  any  
$f \in C({\II}^d)$ and $x \in {\II}^d, $ $Q(f,x)$ 
depends only on the value $f(y)$ at an absolute constant number of points 
 $y$ with $|y_i - x_i| \le \delta 2^{-k_i}, i \in N[d]$;
\begin{equation} \label{ineq:Boundedness}
\|Q_k(f)\|_{C({\II}^d)}  \le C \|\Lambda \|^d \|f\|_{C({\II}^d)} 
\end{equation}
for each $f \in C({\II}^d)$ with a constant $C$ not depending on $k$; and, 
\begin{equation} \label{eq:Reproducing^d}
Q_k(p^*) \ = \ p, \ p \in \Pp_{r-1}^d, 
\end{equation}
where $p^*$ is the restriction of $p$ on ${\II}^d.$ 
The multivariate $Q_k$ is called 
{\it a mixed quasi-interpolant in }$C({\II}^d).$ 

From \eqref{ineq:Boundedness} and \eqref{eq:Reproducing^d} we can see that
\begin{equation} \label{ConvergenceMixedQ_k(f)}
\|f - Q_k(f)\|_{C({\II}^d)}  \to 0 , \ k  \to  \infty.
\end{equation}
(Here and in what follows, $k  \to  \infty$ means that $k_i  \to  \infty$ for $i \in N[d]$).

If $k \in {\Bar{\ZZ}}^d_+$,
we define $T_k := I - Q_k$ for the univariate operator $Q_k$, 
where $I$ is the identity operator. if $k \in {\Bar{\ZZ}}^d_+$, we define the mixed operator
$T_k$ in the manner of the definition \eqref{def:Mixed[Q_k]} by
\begin{equation*}  
T_k:= \prod_{i=1}^d  T_{k_i}.
\end{equation*}
For any $e \subset N[d]$, put 
${\Bar{\ZZ}}^d_+(e):= \{s \in {\Bar{\ZZ}}^d_+: s_i > - 1 , \ i \in e, \ s_i = - 1 , \ i \notin e\}$
(in particular, ${\Bar{\ZZ}}^d_+(\emptyset)= \{(-1,-1,...,-1)\}$ and 
${{\Bar{\ZZ}}}^d_+(N[d])= {\ZZ}^d_+$).
 We have ${\Bar{\ZZ}}^d_+(u)\cap {\Bar{\ZZ}}^d_+(v) = \emptyset$ if $u \ne v$, and 
the following decomposition of ${\Bar{\ZZ}}^d_+$:
\begin{equation*}
{\Bar{\ZZ}}^d_+ =  \bigcup_{e \subset N[d]} {\Bar{\ZZ}}^d_+(e).
\end{equation*}

If $\tau$ is a number such that $0 < \tau \le \min (p,1),$ then for any sequence of 
functions $\{ g_k \}$ there is the inequality
\begin{equation} \label{ineq:L_pNorm}
\left\| \sum   g_k  \right\|_p^{\tau} 
\  \le \ 
 \sum 
 \| g_k \|_p^{\tau}.
\end{equation}

\begin{lemma}  \label{Lemma:IneqNormMixedT_k(f)}
Let $0 < p \le \infty$ and $\tau \le \min (p,1)$. 
Then for any $f \in C(\II^d)$ and $k \in {\Bar{\ZZ}}^d_+(e)$, there holds the inequality 
\begin{equation} \label{ineq:NormMixedT_k(f)}
\left\| T_k(f)\right\|_p  
\ \le \ 
C   \left( \sum_{s\in {\ZZ}^d_+(e), \ s \ge k}
 \left\{ 2^{|s-k|_1/p}\omega_{r}^e(f,2^{-s})_p \right\}^\tau \right)^{1/\tau}
\end{equation}
with some constant $C$ depending at most on $r, \mu, p, d$ and $\|\Lambda\|,$  
whenever the sum in the right-hand side is finite.
\end{lemma}

\begin{proof} 
Notice that ${\Bar{\ZZ}}^d_+(\emptyset)= \{(-1,-1,...,-1)\}$ 
and consequently, the inequality \eqref{ineq:NormMixedT_k(f)} is trivial for 
$e = \emptyset$: $ \|f\|_p \ \le \ C \omega_{r}^\emptyset (f,1)_p \ = \ C \|f\|_p$. 
Consider the case where $e \ne \emptyset$.
 For simplicity we prove the lemma for $d=2$ and $e = \{1,2\} $, 
i.e., ${\Bar{\ZZ}}^d_+(e)= {\ZZ}^2_+$.
This lemma has proven in \cite{Di6}, \cite{Di7} for univariate functions ($d=1$) and even $r$. 
It can be proven for univariate functions and odd $r$ in a completely similar way.
Therefore, 
by \eqref{ineq:w_l><omega_l} there holds the inequality
\begin{equation} \label{ineq:NormT_k(f)} 
\| T_{k_i}(f))\|_p  
\ \ll \ 
\left( \sum_{s_i \ge k_i}
 \left\{ 2^{(s_i-k_i)/p}  \left( 2^{-s_i} \int_{U(2^{-s_i})} \int_{\II(h_i)}
|\Delta_{h_i}^l(f,x)|^p \ d{x_i} \ d{h_i} \right)^{1/p} \right\}^\tau \right)^{1/\tau}, \ \ i=1,2,
\end{equation}
where the norm
$\| T_{k_i}(f))\|_p$ is applied to the univariate function $f$ by considering $f$ as a 
function of  variable $x_i$ with the other variable held fixed.

If $1 \le p < \infty$, we have by \eqref{ineq:NormT_k(f)} applied for $i=1$,
\begin{equation*} 
\begin{aligned}
\| T_{k_1}T_{k_2}(f))\|_p  
\ &  \ll \ 
\left( \int_{\II} \left\{ \sum_{s_1 \ge k_1}
  2^{(s_1-k_1)/p}  \left( 2^{-s_1}  \int_{U(2^{-s_1})} \int_{\II(h_1)}
|\Delta_{h_1}^l((T^1_{k^2_2}f),x)|^p \ d{x_1} \ d{h_1} \right)^{1/p} \right\}^p dx_2 \right)^{1/p}  \\
\ &  \ll \ 
  \sum_{s_1 \ge k_1}
  2^{(s_1-k_1)/p}  \left( 2^{-s_1} \int_{\II}  \int_{U(2^{-s_1})} \int_{\II(h_1)}
|\Delta_{h_1}^l((T^2_{k_2}f),x)|^p \ d{x_1} \ d{h_1} \ dx_2 \right)^{1/p}   \\
\ &  = \ 
  \sum_{s_1 \ge k_1}
  2^{(s_1-k_1)/p}  \left( 2^{-s_1}  \int_{U(2^{-s_1})} \int_{\II(h_1)}\left\{ \int_{\II} 
|(T_{k_2}(\Delta_{h_1}^lf),x)|^p \ dx_2 \right\}\ d{x_1} \ d{h_1} \right)^{1/p}.
\end{aligned}
\end{equation*}
Hence, applying \eqref{ineq:NormT_k(f)} with $i=2$  gives
\begin{equation*} 
\begin{aligned}
\| T_{k_1}T_{k_2}(f))\|_p  
\ &  \ll \ 
   \sum_{s \ge k}
  2^{|s - k|_1)/p}  \left( 2^{-|s|_1}  \int_{U(2^{-s})} \int_{\II^2(h)}
  |\Delta_h^l(f,x)|^p \ dx dh \right)^{1/p}   \\
\ &  \ll \ 
  \sum_{s \ge k}
  2^{|s - k|_1)/p} w_{r}(f,2^{-k})_p    
 \   \ll \ 
  \sum_{s \ge k}
  2^{|s - k|_1)/p} \omega_{r}(f,2^{-k})_p.  
\end{aligned}
\end{equation*}
Thus, the lemma has proven when $1 \le p < \infty$. The cases $0 < p < 1$ and $p = \infty$ can be proven in a 
similar way.
\end{proof}

Let $J_r^d(k):= J^d(k)$ if $r$ is even, and 
$J_r^d(k):= \{s \in {\ZZ}^d: - r < s_i < 2^{k_i+1} + r, i \in N[d] \}$ 
if $r$ is odd. Notice that $J_r^d(k)$ is the set of $s$ for which 
$M^r_{k,s}$ do not vanish identically on  ${\II}^d$. 
Denote by $\Sigma_r^d(k)$ the span of the B-splines 
$M^r_{k,s}, \ s \in J_r^d(k)$.
If $0 < p \le \infty,$ 
for all $k \in {\ZZ}^d_+$ and all $g \in \Sigma_r^d(k)$ such that
\begin{equation} \label{def:StabIneq}
g = \sum_{s \in J_r^d(k)} a_s M^r_{k,s},
\end{equation}
there is the quasi-norm equivalence 
\begin{equation} \label{eq:StabIneq}
 \|g\|_p 
\ \asymp \ 2^{- |k|_1/p}\|\{a_s\}\|_{p,k},
\end{equation}
where
\begin{equation*}
\| \{a_s\} \|_{p,k}
 \ := \
\biggl(\sum_{s \in J_r^d(k)}| a_s |^p \biggl)^{1/p} 
\end{equation*}
with the corresponding change when $p= \infty.$

Let  the mixed operator $q_k$, $k \in {\ZZ}^d_+$, be defined 
in the manner of the definition \eqref{def:Mixed[Q_k]} by
\begin{equation} \label{eq:Def[q_k]}
q_k\ := \ \prod_{i=1}^d \left(Q_{k_i}- Q_{k_i-1}\right). 
\end{equation}
We have 
 \begin{equation} \label{eq:MixedQ_k(2)}
Q_k
\ = \ 
\sum_{k' \le k}q_{k'}. 
\end{equation}
Here and in what follows, for $k,k' \in {\ZZ}^d$ the inequality $k' \le k$ 
means $k'_i \le k_i, \ i \in N[d]$.
From \eqref{eq:MixedQ_k(2)} and \eqref{ConvergenceMixedQ_k(f)} it is easy to see that 
a continuous function  $f$ has the decomposition
\begin{equation} \label{eq:Decomposition}
f \ =  \ \sum_{k \in {\ZZ}^d_+} q_k(f)
\end{equation}
with the convergence in the norm of $C({\II}^d)$.  

From the definition \eqref{eq:Def[q_k]} and the refinement equation for the B-spline $M$, we can represent 
the component functions $q_k(f)$ as 
 \begin{equation} \label{eq:RepresentationMixedq_k(f)}
q_k(f) 
= \ \sum_{s \in J_r^d(k)}c^r_{k,s}(f)M^r_{k,s},
\end{equation}
where $c^r_{k,s}$ are certain coefficient functionals of 
$f,$ which are defined as follows. We first consider the univariate case.
We have
\begin{equation} \label{eq:[q_k]}
q_k(f) 
\ = \ 
\sum_{s \in J(k)}a_{k,s}(f)M_{k,s} - \sum_{s \in J(k-1)}a_{k-1,s}(f)M_{k-1,s}.
\end{equation}
If the order $r$ of the B-spline $M$ is even, 
by using the  refinement equation
\begin{equation} \label{eq:Refinement}
M(x) 
\ = \ 
2^{-r+1} \sum_{j=0}^r \binom{r}{j} M(2x - j + r/2),  
\end{equation}
from \eqref{eq:[q_k]} we obtain 
\begin{equation} \label{eq:Representation[q_k(f)]}
q_k(f) 
\ = \ 
\sum_{s \in J_r(k)}c^r_{k,s}(f)M^r_{k,s},
\end{equation}
where 
\begin{equation}\label{def:c_{k,s}}
c^r_{k,s}(f)
\ := \
a_{k,s}(f) - a^\prime_{k,s}(f), \ k > 0,
\end{equation}
\begin{equation*}
a^\prime_{k,s}(f):= \ 2^{-r+1} \sum_{(m,j) \in C(k,s)}
 \binom{r}{j} a_{k-1,m}(f),  \ k > 0, \ \ 
a^\prime_{0,s}(f):= 0.
\end{equation*}
 and 
\begin{equation*} 
C_r(k,s) := \{(m,j): 2m + j - r/2 = s, \ m \in J(k-1), \  0 \le j \le r\}, \ k > 0, \ \
C_r(0,s) := \{0\}.
\end{equation*} 

If the order $r$ of the B-spline $M$ is odd, 
by using \eqref{eq:Refinement}
from \eqref{eq:[q_k]} we get  \eqref{eq:Representation[q_k(f)]} with
\begin{equation*}
c^r_{k,s}(f)
\ := \
\begin{cases}
0, & \ k = 0, \\
a_{k,s/2}(f), & \ k >0, \ s \ \text{even}, \\
 2^{-r+1} \sum_{(m,j) \in C_r(k,s)}
 \binom{r}{j} a_{k-1,m}(f), & \ k >0, \ s \ \text{odd}, 
\end{cases}
\end{equation*}
 where
\begin{equation*} 
C_r(k,s) := \{(m,j): 4m + 2j - r = s, \ m \in J(k-1), \  0 \le j \le r\}, \ k > 0, \ \
C_r(0,s) := \{0\}.
\end{equation*} 
In the multivariate case, the representation  \eqref{eq:RepresentationMixedq_k(f)} holds true 
with the $c^r_{k,s}$ which are defined in the manner of the definition
\eqref{def:Mixed[a_{k,s}(f)]} by
\begin{equation} \label{def:Mixedc_{k,s}}
c^r_{k,s}(f) 
\ = \   
c^r_{k_1,s_1}((c^r_{k_2,s_2}(... c^r_{k_d,s_d}(f))).
\end{equation}

Let us use the notations: ${\bf 1}:= (1,...,1) \in \RR^d$; $x_+:= ((x_1)_+, ..., (x_d)_+)$ for $x \in \RR^d$; 
${\NN}^d(e):= \{s \in {\ZZ}^d_+: s_i > 0 , \ i \in e, \ s_i = 0 ,  i \notin e\}$ for $e \subset N[d]$ 
(in particular, ${\NN}^d(\emptyset)= \{0\}$ and ${\NN}^d(N[d])= {\NN}^d$).
 We have ${\NN}^d(u)\cap {\NN}^d(v) = \emptyset$ if $u \ne v$, 
 and the following decomposition of ${\ZZ}^d_+$:
\begin{equation*}
{\ZZ}^d_+ =  \bigcup_{e \subset N[d]} {\NN}^d(e).
\end{equation*}

\begin{lemma} \label{Lemma:IneqMixedq_k(f))} 
Let $0 < p \le \infty$ and $\tau \le \min (p,1)$. 
Then for any $f \in C(\II^d)$ and $k \in {\NN}^d(e)$, there holds the inequality  
\begin{equation*}  
\| q_k(f))\|_p  
\ \le \ 
C  \sum_{v \supset e} \left( \sum_{s \in {\ZZ}^d_+(v), \ s \ge k}
 \left\{ 2^{|s-k|_1/p}\omega_{r}^v(f,2^{-s})_p \right\}^\tau \right)^{1/\tau}
\end{equation*}
with some constant $C$ depending at most on $r, \mu, p, d$ and $\|\Lambda\|,$  
whenever the sum in the right-hand side is finite.
\end{lemma}
\begin{proof}
From the equality 
\begin{equation*}
q_k = \prod_{i=1}^d \left(T^i_{k_i-1} - T^i_{k_i}\right), 
\end{equation*}
it follows that 
 \begin{equation*}  
 q_k  
\  = \ 
 \sum_{u \subset N[d]} (-1)^{|u|} \prod_{i \in u} T^i_{k_i} \prod_{i \notin u}T^i_{k_i - 1} 
\  = \ 
\sum_{u \subset N[d]} (-1)^{|u|} T_{k^u},
\end{equation*}
where $k^u$ is defined by $k^u_i = k_i$ if $i \in u$, and 
$k^u_i = k_i - 1$ if $i \notin u$. 
 Hence,
 \begin{equation} \label{ineq:Mixedq_k (f)} 
\|q_k (f)\|_p 
\ \ll \ 
\sum_{u \subset N[d]} \| T_{k^u}(f)\|_p.
\end{equation}
Notice that $k^u \in {\Bar{\ZZ}}^d_+(v)$ for some $v \supset e$, 
and $0 \ \le \ k - k^u_+ \ \le k - k^u \ \le \ {\bf 1}$. Moreover, for $s \in {\ZZ}^d_+(v)$, $s \ge k^u$ 
if only if $s \ge k^u_+$. 
Hence, by Lemma \ref{Lemma:IneqNormMixedT_k(f)} and properties of mixed modulus smoothness we have
\begin{equation*}  
\begin{aligned}
\| T_{k^u}(f))\|_p  
\ &  \ll \   
  \left( \sum_{s \in {\ZZ}^d_+(v), s \ge k^u} 
 \left\{ 2^{|s-k^u|_1/p}\omega_{r}^v(f,2^{-s})_p \right\}^\tau \right)^{1/\tau} \\
\ &  \ll \   
  \left( \sum_{s \in {\ZZ}^d_+(v), s \ge k^u_+} 
 \left\{ 2^{|s-k^u_+|_1/p}\omega_{r}^v(f,2^{-s})_p \right\}^\tau \right)^{1/\tau} \\
\ &  = \   
 \left( \sum_{s' \in {\ZZ}^d_+(v), s' \ge k}
 \left\{ 2^{|s' - k|_1/p}\omega_{r}^v(f,2^{- (s' - (k - k^u_+))})_p \right\}^\tau \right)^{1/\tau} \\
 \ &  \ll \   
  \left( \sum_{s \in {\ZZ}^d_+(v), s \ge k}
 \left\{ 2^{|s - k|_1/p}\omega_{r}^v(f,2^{- s})_p \right\}^\tau \right)^{1/\tau}.
\end{aligned}
\end{equation*}
The last inequality together with \eqref{ineq:Mixedq_k (f)} proves the lemma.
\end{proof}

\begin{lemma} \label{Lemma:IneqMixedomega_{r}^e}  
Let $0 < p \le \infty, \ 0 < \tau \le \min(p,1), \ \delta = \min (r, r - 1 + 1/p).$ 
Then for any $f \in C(\II^d)$ and $k \in {\ZZ}^d_+(e)$, there holds the inequality 
\begin{equation*} 
\omega_{r}^e(f,2^{-k})_p 
\ \le \ 
C   \left( \sum_{s\in {\ZZ}^d_+}
 \left\{ 2^{-\delta |(k-s)_+|_1} \| q_s(f))\|_p \right\}^\tau \right)^{1/\tau}
\end{equation*}
with some constant $C$ depending at most on $r, \mu, p, d$ and $\|\Lambda\|,$  
whenever the sum in the right-hand side is finite.
\end{lemma}

\begin{proof}
For simplicity we prove the lemma for $e= N[d]$, i.e., ${\ZZ}^d_+(e)= {\ZZ}^d_+$. 
Let $f \in C({\II}^d)$ and $k \in {\ZZ}^d_+$. From \eqref{eq:Decomposition}
and  \eqref{ineq:L_pNorm} we obtain
\begin{equation} \label{ineq:NormDelta_h^{r}(f)}  
\|\Delta_h^{r}(f)\|_p 
\ \le \ 
C   \left( \sum_{s \in {\ZZ}^d_+}
  \| \Delta_h^{r}(q_s(f))\|_p ^\tau \right)^{1/\tau}.
\end{equation}
Further, by \eqref{eq:RepresentationMixedq_k(f)} we get
\begin{equation*} 
\Delta_h^{r}(q_k(f)) 
= \ \sum_{j \in J^d(s)}c^r_{s,j}(f)\Delta_h^{r}(M^r_{s,j}).
\end{equation*}
Notice that for any $x$, the number of non-zero B-spines in \eqref{eq:RepresentationMixedq_k(f)} 
is an absolute constant depending on $r,d$ only. Thus, we have 
\begin{equation} \label{ineq:|Delta_h^{r}(q_s(f))|^p}
|\Delta_h^{r}(q_s(f),x)|^p 
\ \ll \ 
\sum_{j \in J^d(s)}|c^r_{s,j}(f)|^p|\Delta_h^{r}(M^r_{s,j},x)|^p, \ \ x \in {\II}^d.
\end{equation} 
From properties of  the B-spline $M$ it is easy to prove the following estimate
 \begin{equation*} 
\int_{{\II}^d(h)}|\Delta_h^{r}(M^r_{s,j},x)|^p dx 
\ \ll \
 2^{- |s|_1 - \delta p |(- \log|h| - s)_+|_1},
\end{equation*}
where we used the abbreviation $\log|h|:= (\log |h_1|, ..., \log |h_d|)$.
Hence, by \eqref{ineq:|Delta_h^{r}(q_s(f))|^p} we obtain
 \begin{equation*} 
 \begin{aligned}
\|\Delta_h^{r}(q_s(f))\|_p 
\ & \ll \
 2^{- \delta |(- \log|h| - s)_+|_1}
2^{- |s|_1} \left(\sum_{j \in J^d(s)}|c^r_{s,j}(f)|^p \right)^{1/p} \\
\ & \ll \
 2^{- \delta |(- \log|h| - s)_+|_1}\|q_s(f)\|_p.
\end{aligned}
\end{equation*}
By \eqref{ineq:NormDelta_h^{r}(f)} we have
\begin{equation*}  
\|\Delta_h^{r}(f)\|_p  
\ \ll \ 
  \left( \sum_{s \in {\ZZ}^d_+}
 \left\{ 2^{- \delta |(- \log|h|-s)_+|_1} \| q_s(f))\|_p \right\}^\tau \right)^{1/\tau}.
\end{equation*}
From the last inequality we prove the lemma.
\end{proof}

For functions $f$ on ${\II}^d$, we introduce the following quasi-norms:
\begin{equation*} 
\begin{aligned}
B_2(f)
 \ & := \ 
\| \{ q_k(f) \}\|_{b_\theta^\alpha(L_p)}; \\
B_3(f)
\ & := \ 
\biggl(\sum_{k=0}^\infty 
\bigl( 2^{(\alpha - d/p)k}\|\{c^r_{k,s}(f)\}\|_{p,k} \bigl)^\theta \biggl)^{1/\theta}.
\end{aligned}
\end{equation*}

We will need the following discrete Hardy inequality. 
Let $\{a_k\}_{k \in {\ZZ}^d_+}$ and $\{b_k\}_{k \in {\ZZ}^d_+}$ be two positive 
sequences and let for some $M > 0, \ \tau >0$ 
\begin{equation} \label{HardyIneq(1)}
b_k
 \  \le \ M \left(  \sum_{s \in {\ZZ}^d_+} 
\left(2^{\delta |(k-s)_+|_1}a_s\right)^\tau \right)^{1/\tau}.
\end{equation}
Then for any $0 < \beta < \delta, \ \theta >0,$ 
\begin{equation} \label{HardyIneq(2)}
\ \|\{ b_k\}\|_{b_\theta^\beta}
\ \le \ 
 C M \|\{ a_k\}\|_{b_\theta^\beta}
\end{equation}
with $C = C(\beta, \theta, d)$. For a proof of this inequality for the univariate case 
see, e.g, \cite{DL}. In the general case it can be proven by induction based on the 
univariate case.

\begin{theorem} \label{Theorem:Representation}
Let $\ 0 < p, \theta \le \infty$ and $1/p < \alpha < r$. 
Then the hold the following assertions. 
\begin{itemize}
\item[(i)] A function $f \in \MB$ can be represented by the mixed B-spline series 
\begin{equation} \label{eq:B-splineRepresentation}
f \ = \sum_{k \in {\ZZ}^d_+} \ q_k(f) = 
\sum_{k \in {\ZZ}^d_+} \sum_{s \in J_r^d(k)} c^r_{k,s}(f)M^r_{k,s}, 
\end{equation}  
satisfying 
the convergence condition
\begin{equation*} 
B_2(f) \ \asymp \ B_3(f) \ \ll \ B(f),
\end{equation*}
where the coefficient functionals $c^r_{k,s}(f)$ are explicitly constructed by formula 
\eqref{def:c_{k,s}}--\eqref{def:Mixedc_{k,s}} as
linear combinations of an absolute constant number of values of 
$f$ which does not depend on neither $k,s$ nor $f$.  
\item[(ii)]  If in addition, $ \alpha < \min(r, r - 1 + 1/p)$,
then a continuous function $f$ on ${\II}^d$ belongs to 
the Besov space $\MB$ if and only if 
$f$ can be represented by the series \eqref{eq:B-splineRepresentation}.
Moreover, the Besov quasi-norm $B(f)$ is equivalent to one of the quasi-norms $B_2(f)$ and $B_3(f)$.
\end{itemize} 
\end{theorem}

\begin{proof} Since by  \eqref{eq:StabIneq} 
the quasi-norms $B_2(f)$ and $B_3(f)$ are equivalent,
it is enough to prove (i) and (ii) for $B_3(f)$. Fix a number $0 < \tau \le \min (p,1).$ 

\noindent
{\it Assertion} (i): For $k \in {\ZZ}^d_+$, put 
\begin{equation*}
b_k := 2^{|k|_1/p} \|q_k (f)\|_p, \ \
 a_k 
\ := \ 
\left( \sum_{v \supset e}
 \left\{ 2^{|k|_1/p}\omega_{r}^v(f,2^{-k})_p \right\}^\tau \right)^{1/\tau}
\end{equation*}
 if $k \in {\NN}^d(e)$.
 By Lemma \ref{Lemma:IneqMixedq_k(f))} we have for $k \in {\ZZ}^d_+$,
\begin{equation*}
b_k 
\ \le \ 
C   \left( \sum_{s \ge k}^\infty
 a_s ^{\tau} \right)^{1/\tau}.
\end{equation*}
Then applying the mixed discrete Hardy inequality \eqref{HardyIneq(1)}--\eqref{HardyIneq(2)} 
with $\beta = \alpha - 1/p$, gives 
\begin{equation*} 
B_3(f) \ = \ \|\{ b_k\}\|_{b_\theta^\beta}
\ \le \ 
 C \|\{ a_k\}\|_{b_\theta^\beta} \ \asymp \ B_1(f) \ \asymp \ B(f).
\end{equation*}

\noindent
{\it Assertion} (ii): Let in addition, $ \alpha < \min (r, r - 1 + 1/p)$.
For $k \in {\ZZ}^d_+$, put 
\begin{equation*}
b_k := \left( \sum_{v \supset e}
 \left\{ \omega_{r}^v(f,2^{-k})_p \right\}^\tau \right)^{1/\tau} \ \
 a_k 
\ := \ 
\|q_k (f)\|_p
\end{equation*}
 if $k \in {\NN}^d(e)$. 
By Lemma \ref{Lemma:IneqMixedomega_{r}^e} we have for $k \in {\ZZ}^d_+$
\begin{equation*} 
b_k
 \  \le \ C \left(  \sum_{s \in {\ZZ}^d_+} 
\left(2^{\delta |(k-s)_+|_1}a_s\right)^\tau \right)^{1/\tau},
\end{equation*}
where $\delta = \min (r, r - 1 + 1/p)$.
Then applying the mixed discrete Hardy inequality \eqref{HardyIneq(1)}--\eqref{HardyIneq(2)}  
with $\beta = \alpha$, gives 
\begin{equation*} 
B(f) \ \asymp \ B_1(f) \ \asymp \ \|\{ b_k\}\|_{b_\theta^\beta}
\ \le \ 
 C \|\{ a_k\}\|_{b_\theta^\beta} \ = \ B_3(f).
\end{equation*}
The assertion (ii) is proven.
\end{proof}

\noindent
{\bf Remark}\ From \eqref{def:c_{k,s}}--\eqref{def:Mixedc_{k,s}} we can see that if $r$ is even,
for each pair $k,s$ the coefficient $c^r_{k,s}(f)$ is a linear combination of the values 
$f(2^{-k}(s - j)),$  and $f(2^{-k + {\bf 1}}(s^\prime - j)), \ j \in P^d(\mu), \ s^\prime \in C_r(k,s)$.
The number of these values does not exceed the fixed number $(2 \mu + 1)^d((r + 1)^d + 1)$. 
If $r$ is odd, we can say similarly about the coefficient $c^r_{k,s}(f)$.

\section{Sampling recovery} 
\label{Sampling recovery}

\medskip
Recall that  the linear operator $R_m, m \in {\ZZ}_+$, is defined for functions on ${\II}^d$ 
in \eqref{def:R_m} as follows.
\begin{equation} \label{eq:R_m}
R_m(f) 
\ = \ 
\sum_{k \in \Delta(m)} q_k(f)
\ = \ 
\sum_{k \in \Delta(m)} \sum_{s \in J_r^d(k)} c^r_{k,s}(f)M^r_{k,s}. 
\end{equation} 

\begin{lemma} \label{Lemma:R_m=L_n}
 For functions $f$ on ${\II}^d$, $R_m$ defines a linear sampling algorithm 
of the form \eqref{def:L_n(f)} on the grid $G^d(m)$. More precisely, 
\begin{equation*} 
R_m(f) 
\ = \ 
L_n(f) 
\ = \ 
\sum_{(k,s) \in G^d(m)} f(2^{-k}j) \psi_{k,j}, 
\end{equation*} 
where 
\begin{equation} \label{def:[n]}
n \ := \ |G^d(m)|
\ = \ 
\sum_{k \in \Delta(m)} |I^d(k)|
\ \asymp \ 
2^m  m^{d-1}; 
\end{equation} 
$\psi_{k,j}$ are explicitly constructed as linear combinations of at most 
$(4\mu + r + 5)^d$ B-splines $M^r_{k,s} \in  M_r^d(m)$
for even $r$, and $(12\mu + 2r + 13)^d$ B-splines  $M^r_{k,s} \in  M_r^d(m)$ for odd $r$;
$M_r^d(m) := \{M^r_{k',s'}: k' \in \Delta(m), s' \in J_r^d(k')\}$. 
\end{lemma}

\begin{proof}
Let us prove the lemma for even $r$. For odd $r$ it can be proven in a similar way. 
 For univariate functions the coefficient functionals 
$a_{k,s}(f)$ can be rewritten as 
\begin{equation*} 
a_{k,s}(f)
\  = \   
\sum_{|s - j| \le \mu} \lambda (s - j) f_k(2^{-k}j)
\  = \   
\sum_{j \in P(k,s)} \lambda_{k,s} (j) f(2^{-k}j),
\end{equation*}
where $\lambda_{k,s} (j) := \lambda(s - j)$ and
$P(k,s) = P_s(\mu) := \{j \in \{0,2^k\}: s-j \in P(\mu)\}$ for $\mu \le s \le 2^k - \mu$;  
$\lambda_{k,s} (j)$ is  a linear combination of no more than
$\max (r, 2\mu +1)\le 2\mu + 2$ coefficients $\lambda (j), j \in P(\mu)$, for
$s < \mu$ or $s > 2^k - \mu$ and
\begin{equation*}
P(k,s)
\ \subset \
\begin{cases} 
P_s(\mu) \cup \{0,r - 1\}, \ & s < \mu, \\
P_s(\mu)\cup \{2^k - r + 1, 2^k\}, \ & s > 2^k - \mu.
\end{cases}
\end{equation*}
Further, for univariate functions we have
\begin{equation*} 
\begin{aligned}
c^r_{k,s}(f)
\ & = \   
\sum_{j \in P(k,s)} \lambda_{k,s} (j) f(2^{-k}j) - 
2^{-r+1} \sum_{(m,\nu) \in C_r(k,s)} \binom{r}{\nu}
\sum_{j \in P(k-1,m)} \lambda_{k-1,m} (j) f(2^{-k}(2j)) \\
\ & = \   
\sum_{j \in G(k,s)} \lambda_{k,s} (j) f(2^{-k}j),
\end{aligned}
\end{equation*}
where $G(k,s):= P(k,s) \cup \{2j: j \in P(k-1,m), \ (m,\nu)\in C(k,s) \}$. 
If $j \in P(k,s)$, we have $|j - s|\le \max(r, 2\mu + 1) \le 2\mu + 2$. 
If $j \in P(k-1,m), \ (m,\nu)\in C(k,s)$, we have 
$|2j - s| = |2j - 2m - \nu + r/2| \le 2|j - m| +|\nu - r/2| 
\le 2\max(r, 2\mu + 1) + r + 1 \le 4\mu + r + 5 =: {\bar \mu}$. Therefore, 
$G(k,s) \subset P_s({\bar \mu})$, and we can rewrite the coefficient functionals $c^r_{k,s}(f)$
in the form 
\begin{equation*} 
c^r_{k,s}(f)
\  = \   
\sum_{j - s \in P({\bar \mu})} \lambda_{k,s} (j) f(2^{-k}j)
\end{equation*}
with zero coefficients $\lambda_{k,s} (j)$ for $j \notin G(k,s)$. 
Therefore, we have  
\begin{equation*}
\begin{aligned} 
q_k(f)
\ & = \
\sum_{s \in J_r(k)} c^r_{k,s}(f) M^r_{k,s} 
\  = \   
\sum_{s \in J_r(k)} \sum_{j - s \in P({\bar \mu})} \lambda_{k,s} (j) f(2^{-k}j) M^r_{k,s} \\
\ & = \   
\sum_{j \in I(k)} f(2^{-k}j) \sum_{s-j \in P({\bar \mu})} \gamma_{k,j} (s)  M^r_{k,s}
\end{aligned}
\end{equation*}
for certain coefficients $\gamma_{k,j}(s)$. Thus,  the univariate $q_k(f)$ is of the form
\begin{equation*} 
q_k(f)
\ = \
\sum_{j \in I(k)} f(2^{-k}j) \psi_{k,j},
\end{equation*}
where
\begin{equation*}
 \psi_{k,j}
\ := \
\sum_{s-j \in P({\bar \mu})} \gamma_{k,j} (s)  M_{k,s},
\end{equation*}
are a linear combination of no more than the absolute number $4\mu + r + 5$ of B-splines $ M^r_{k,s}$, 
and the size $|I(k)|$ is $2^k$.
Hence, the multivariate $q_k(f)$ is of the form
\begin{equation*}
q_k(f)
\ = \
\sum_{j \in I^d(k)} f(2^{-k}j) \psi_{k,j},
\end{equation*}
where
\begin{equation*}
 \psi_{k,j}
\ := \
\prod_{i=1}^d  \psi_{k_i,j_i}
\end{equation*}
are a linear combination of no more than the absolute number 
$(4\mu + r + 5)^d$ of B-splines $ M^r_{k,s} \in M_r^d(m)$, 
and the size $|I^d(k)|$ is $2^{|k|_1}$. From \eqref{eq:R_m} 
we can see that $R_m(f)$ is of the form \eqref{def:L_n(f)} with $n$ as in \eqref{def:[n]}.
\end{proof}

\begin{theorem} \label{Theorem:UpperBoundE(m),alpha>1/p}
Let $\ 0 < p, q, \theta \le \infty$ and $1/p < \alpha < r$. Then we have the following upper bound of $E(m)$. 
\begin{itemize}
\item[{\rm (i)}] For $p \ge q$,
\begin{equation*} 
E(m)
 \ \ll \ 
\begin{cases}
2^{- \alpha m}, \ & \theta \le \min(q,1), \\
2^{- \alpha m} m^{(d-1)(1/q - 1/\theta)}, \ & \theta > \min(q,1), \ q \le 1, \\
2^{- \alpha m} m^{(d-1)(1 - 1/\theta)}, \ & \theta > \min(q,1), \ q > 1.
\end{cases}
\end{equation*}
\item[{\rm (ii)}] For $p < q$, 
\begin{equation*} 
E(m)
 \ \ll \ 
\begin{cases}
2^{- (\alpha - 1/p + 1/q) m} m^{(d-1)(1/q - 1/\theta)_+}, \ & q < \infty, \\
2^{- (\alpha - 1/p) m} m^{(d-1)(1 - 1/\theta)_+}, \ &  q = \infty.
\end{cases}
\end{equation*}
\end{itemize} 
\end{theorem}

\begin{proof} 

\noindent
{\it Case} (i): $p \ge q$. 
For an arbitrary $f \in B^\alpha_{p,\theta}$, by the representation \eqref{eq:B-splineRepresentation}
and \eqref{ineq:L_pNorm} we have
\begin{equation*}
\|f - R_m(f)\|_q^{\tau} 
\  \ll \ 
 \sum_{|k|_1 > m} \| q_k(f)\|_q^{\tau}
\end{equation*}
with any $\tau \le \min(q,1)$. Therefore, if $\theta \le \min(q,1)$, then
by the inequality $\| q_k(f)\|_q \le \| q_k(f)\|_p$ we get 
\begin{equation*} 
\begin{aligned}
\|f - R_m(f)\|_q  
\  & \ll \ 
 \left(\sum_{|k|_1 > m} \| q_k(f)\|_q^{\theta} \right)^{1/\theta} \\
\  & \le \ 
2^{-\alpha m} \left(\sum_{|k|_1 > m} \{2^{\alpha |k|_1}\| q_k(f)\|_p\}^{\theta} \right)^{1/\theta} \\
 \  & \ll \ 
 2^{-\alpha m} B_2(f)
 \ \ll \ 
 2^{-\alpha m}.
\end{aligned}
\end{equation*}
Further, if $\theta > \min(q,1)$, then 
\begin{equation*} 
\|f - R_m(f)\|_q^{q^*}
\  \ll \ 
\sum_{|k|_1 > m} \| q_k(f)\|_q^{q^*}
 \  \ll \ 
 \sum_{|k|_1 > m} \{2^{\alpha |k|_1}\| q_k(f)\|_q\}^{q^*}\{2^{- \alpha |k|_1}\}^{q^*},
\end{equation*}
where $q^* = \min(q,1)$. Putting $\nu = \theta/q^*$, by H\"older's inequality and 
the inequality $\| q_k(f)\|_q \le \| q_k(f)\|_p$ we obtain
\begin{equation} \label{ineq:E^{q^*}(m)}
\begin{aligned}
\|f - R_m(f)\|_q^{q^*}
\ & \ll \ 
\left( \sum_{|k|_1 > m} \{2^{\alpha |k|_1}\| q_k(f)\|_q\}^{q^* \nu}\right)^{1/\nu} 
\left(\sum_{|k|_1 > m}\{2^{- \alpha |k|_1}\}^{q^* \nu'}\right)^{1/\nu'} \\
\ & \ll \ 
\{B_2(f)\}^{q^*} \{2^{- \alpha m} m^{(d-1)(1/q^* - 1/\theta)}\}^{q^*} 
\  \ll \ 
\{2^{- \alpha m} m^{(d-1)(1/q^* - 1/\theta)}\}^{q^*}. 
\end{aligned}
\end{equation}
This proves the Case (i).

\noindent
{\it Case} (ii): $p < q$.  We first assume  $q < \infty$.
For an arbitrary $f \in B^\alpha_{p,\theta}$, by the representation \eqref{eq:B-splineRepresentation} 
and Lemma \ref{Lemma:IneqL_qNorm<L_pNorm} we have
\begin{equation*} 
\|f - R_m(f)\|_q^q 
\  \ll \ 
 \sum_{|k|_1 > m}\{2^{(1/p - 1/q) |k|_1} \| q_k(f)\|_p\}^q.
\end{equation*}
 Therefore, if $\theta \le q$, then
\begin{equation*}
\begin{aligned} 
\|f - R_m(f)\|_q 
\ & \ll \ 
 \left(\sum_{|k|_1 > m} \{2^{(1/p - 1/q) |k|_1} \| q_k(f)\|_p\}^{\theta} \right)^{1/\theta} \\
 \ & \ll \ 
 2^{- (\alpha - 1/p + 1/q) m} B_2(f)
 \ \ll \ 
 2^{- (\alpha - 1/p + 1/q) m}.
 \end{aligned}
\end{equation*}
Further, if $\theta > q$, then
\begin{equation*}
\begin{aligned} 
\|f - R_m(f)\|_q^q 
\ & \ll \ 
\sum_{|k|_1 > m} \{2^{(1/p - 1/q) |k|_1} \| q_k(f)\|_p\}^q \\
 \ & \ll \ 
 \sum_{|k|_1 > m} \{2^{\alpha |k|_1}\| q_k(f)\|_p\}^q\{2^{- (\alpha - 1/p + 1/q)|k|_1}\}^q.
 \end{aligned}
\end{equation*}
Hence, similarly to \eqref{ineq:E^{q^*}(m)}, we get
\begin{equation*} 
E^q(m)
\  \ll \ 
\{2^{- (\alpha - 1/p + 1/q) m} m^{(d-1)(1/q - 1/\theta)}\}^q. 
\end{equation*}

When $q = \infty$, the Case (ii) can be proven analogously by use the inequality
\begin{equation*} 
\|f - R_m(f)\|_\infty 
\  \ll \ 
 \sum_{|k|_1 > m}2^{|k|_1/p} \| q_k(f)\|_p.
\end{equation*}
\end{proof}

The following theorem for the case $\alpha = 1/p$ can be proven by use of 
Lemmas \ref{Lemma:IneqMixedq_k(f))} and \ref{Lemma:IneqL_qNorm<L_pNorm} and the inequality \eqref{ineq:L_pNorm}.

\begin{theorem} \label{Theorem:UpperBoundE(m),alpha=1/p}
Let $\ 0 < p, q < \infty$, $\ 0 < \theta \le \min (p,1)$ 
and $\alpha = 1/p < r$. Then we have the following upper bound of $E(m)$. 
\begin{itemize}
\item[{\rm (i)}] For $p \ge q$,
\begin{equation*} 
E(m) 
\ \ll \ 
\begin{cases}
2^{- m/p} m^{(d-1)}, \ &  \ p \ge 1, \\ 
2^{- m/p} m^{(d-1)/p}, \ & \ p < 1. 
\end{cases}
\end{equation*}
\item[{\rm (ii)}] For $p < q$,
\begin{equation*} 
E(m) 
\ \ll \ 
2^{- m/q} m^{(d-1)/q}. 
\end{equation*}
\end{itemize} 
\end{theorem}

The following two theorems are a direct corollary of Lemma \ref{Lemma:R_m=L_n} 
and Theorems \ref{Theorem:UpperBoundE(m),alpha>1/p} and \ref{Theorem:UpperBoundE(m),alpha=1/p}.
 
\begin{theorem} \label{Theorem:UpperBoundr_n,alpha>1/p}
Let $\ 0 < p, q, \theta \le \infty$ and $1/p < \alpha < r$. 
If ${\bar m}$ is the largest integer of  $m$ such that 
\begin{equation*} 
2^m m^{d-1} 
 \asymp \ 
\sum_{k \in \Delta(m)} |I(k)|
\ \le \ 
n, 
\end{equation*} 
then we have the following upper bound of $r_n$ and $E({\bar m})$. 
\begin{itemize}
\item[{\rm (i)}] For $p \ge q$,
\begin{equation*} 
r_n
\ \ll \
E({\bar m})
\ \ll \ 
\begin{cases}
(n^{-1} \log^{d-1}n)^\alpha, \ & \theta \le \min(q,1), \\
(n^{-1} \log^{d-1}n)^\alpha (\log^{d-1}n)^{ 1/q - 1/\theta}, \ & \theta > \min(q,1), \ q \le 1, \\
(n^{-1} \log^{d-1}n)^\alpha (\log^{d-1}n)^{ 1 - 1/\theta}, \ & \theta > \min(q,1), \ q > 1.
\end{cases}
\end{equation*}
\item[{\rm (ii)}] For $p < q$, 
\begin{equation*} 
r_n
\ \ll \
E({\bar m})
 \ \ll \  
\begin{cases}
(n^{-1} \log^{d-1}n)^{\alpha - 1/p + 1/q}(\log^{d-1}n)^{(1/q - 1/\theta)_+}, \ & q < \infty, \\
(n^{-1} \log^{d-1}n)^{\alpha - 1/p}(\log^{d-1}n)^{(1 - 1/\theta)_+}, \ &  q = \infty.
\end{cases}
\end{equation*}
\end{itemize} 
\end{theorem}

\begin{theorem}
Let $\ 0 < p, q < \infty$, $\ 0 < \theta \le \min (p,1)$ 
and $\alpha = 1/p < r$. If ${\bar m}$ is the largest integer of  $m$ such that 
\begin{equation*} 
2^m m^{d-1} 
 \asymp \ 
\sum_{k \in \Delta(m)} |I(k)|
\ \le \ 
n, 
\end{equation*} 
then we have the following upper bound of $r_n$ and $E({\bar m})$.   
\begin{itemize}
\item[{\rm (i)}] For $p \ge q$,
\begin{equation*} 
r_n
\ \ll \
E({\bar m})
 \ \ll \ 
\begin{cases}
(n^{-1} \log^{d-1}n)^{1/p} \log^{d-1}n, \ &  \ p \ge 1, \\ 
(n^{-1} \log^{d-1}n)^{1/p}(\log^{d-1}n)^{ 1/p}, \ & \ p < 1. 
\end{cases}
\end{equation*}
\item[{\rm (ii)}] For $p < q$,
\begin{equation*} 
r_n
\ \ll \
E({\bar m})
 \ \ll \ 
(n^{-1} \log^{d-1}n)^{1/q}(\log^{d-1}n)^{ 1/q}. 
\end{equation*}
\end{itemize} 
\end{theorem}

From Theorem \ref{Theorem:UpperBoundr_n,alpha>1/p} and Lemma \ref{Lemma: Asymp[lambda_n]}
we obtain the following theorem.

\begin{theorem}
Let $ 1 \le p, q \le \infty$, $ 0 < \theta \le \infty$ 
and $1/p < \alpha  < r$. If ${\bar m}$ is the largest integer of  $m$ such that 
\begin{equation*} 
2^m m^{d-1} 
 \asymp \ 
\sum_{k \in \Delta(m)} |I(k)|
\ \le \ 
n, 
\end{equation*} 
then we have the following asymptotic order of $r_n$ and $E({\bar m})$.     
\begin{itemize}
\item[{\rm (i)}] For $p \ge q$ and $\theta \le 1$,
\begin{equation*} 
E({\bar m})
 \ \asymp \ 
r_n
 \ \asymp \ 
 (n^{-1} \log^{d-1}n)^\alpha, \ 
\begin{cases}
 2 \le q < p < \infty, \\
1 < p = q \le \infty.
\end{cases}
\end{equation*}
\item[{\rm (ii)}] For $1 < p < q < \infty$, 
\begin{equation*} 
E({\bar m})
 \ \asymp \ 
r_n
 \ \asymp \ 
(n^{-1} \log^{d-1}n)^{\alpha - 1/p + 1/q}(\log^{d-1}n)^{(1/q - 1/\theta)_+}, \
 \begin{cases}
2 \le p, \ 2 \le \theta \le q,  \\
q \le 2.
\end{cases}
 \end{equation*}
\end{itemize} 
\end{theorem}

\section{Interpolant representations and sampling recovery}
\label{Faber-Schauder} 
 
We first consider a piecewise constant interpolant representation. 
Let $\chi_{[0,1)}$ and $\chi_{[0,1]}$ be the characteristic functions of the half opened and closed 
intervals $[0,1)$ and $[0,1]$, respectively. For $k \in {\ZZ}_+$ and $s = 0,1,..., 2^k - 1$,
we define the system of functions $N_{k,s}$ on $\II$ by  
\begin{equation*} 
N_{k,s}(x)  
\ := \
\begin{cases}
 \chi_{[0,1)}(2^kx - s),  & \, 0 \le s < 2^k - 1, \\  
\chi_{[0,1]}(2^kx - s),  & \,  s = 2^k - 1. 
\end{cases}
\end{equation*}
(In particular, $N_{0,0} = \chi_{[0,1]}$). Obviously, we have for $k > 0$ and 
$s = 0,1,..., 2^k - 1$, 
\begin{equation*} 
N_{k-1,s}
\ = \  
N_{k,2s} + N_{k,2s+1}.
\end{equation*}

We let the operator $\Pi_k$ be defined for functions $f$ on $\II$, for $k \in {\ZZ}_+$, by
\begin{equation*} 
\Pi_k(f) 
\ := \  
 \sum_{s=0}^{2^k-1} f(2^{-k}s) N_{k,s}, \ \text{and} \ \ \Pi_{-1}(f) = 0.
\end{equation*}
Clearly, the linear operator $\Pi_k$ is bounded in $L_\infty(\II)$, reproduces constant functions  
and for any continuous function $f$,
\begin{equation*} 
\|f - \Pi_k(f)\|_\infty 
\ \le \
\omega_1(f, 2^{-k})_\infty,  
 \end{equation*}
and consequently, $\|f - \Pi_k(f)\|_\infty \to 0$, when $k \to \infty$.
Moreover, for any $x \in \II$,  $\Pi_k(f,x)= f(2^{-k}s)$ if $x$ is 
in either the interval $[2^{-k}s,2^{-k}(s + 1))$ for $s=0,1...,2^{-k} - 2$ or the interval 
$[2^{-k}s,1]$ for $s=2^{-k} - 1$, i.e., $\Pi_k$ possesses a local property.
In particular, $\Pi_k(f)$ interpolates $f$ at the points $2^{-k}s$, $s = \{0,1,..., 2^k - 1\}$, that is,
\begin{equation} \label{eq: Interpolation[Pi_k]}
\Pi_k(f,2^{-k}s) 
\ = \  
 f(2^{-k}s), \ s = 0,1,..., 2^k - 1.
\end{equation}
Further, we define for $k \in {\ZZ}_+$,
\begin{equation*} 
\pi_k(f) 
\ := \  
\Pi_k(f) - \Pi_{k-1}(f).
\end{equation*}
From the definition it is easy to check that 
\begin{equation*} 
\pi_k(f)  \ =  
 \sum_{s \in Z_1(k)} \lambda^1_{k,s}(f)\varphi^1_{k,s}, 
\end{equation*}
where $Z_1(0) := \{0\}$, $Z_1(k) := \{0,1,..., 2^{k-1} - 1\}$ for $k > 0$,
\begin{equation*} 
\varphi^1_{k,s}(x):=  N_{k,2s+1}(x),   \, k > 0, \ 
\text{and} \ \varphi^1_{0,0}(x) \ := \ N_{0,0}(x),
\end{equation*}
and
\begin{equation*} 
\lambda^1_{k,s}(f)  
\ := \
 \Delta^1_{2^{-k}}(f,2^{-k+1}s),   \, k > 0,
\ \text{and} \ \lambda^1_{0,0}(f)  
\ := \
 f(0). 
\end{equation*}

We now can see that every $f \in C(\II)$ is represented by the series
\begin{equation} \label{eq:Representation[r=1]}
f 
\ = \  
 \sum_{k \in {\ZZ}_+} \ \sum_{s \in Z_1(k)}\lambda^1_{k,s}(f) \varphi^1_{k,s},
\end{equation}
converging in the norm of $L_\infty(\II)$. 

Next, let us revisit  the univariate piecewise linear (nodal) quasi-interpolant
 for functions on $\RR$ defined in \eqref{def:Q,r=2}
with $M(x)\ = \ (1 - |x|)_+ \ (r=2)$. 
Consider the generated from it by the formula \eqref{eq:Q_k} quasi-interpolant for functions on $\II$
\begin{equation} \label{eq:NodalQ_k}
Q_k(f,x)= \ \sum_{s \in J(k)} f(2^{-k}s) M_{k,s}(x), 
\end{equation}
and the related quasi-interpolant representation
\begin{equation} \label{eq:B-splineRepresentation[r=2]}
f \ = \sum_{k \in {\ZZ}_+} \ q_k(f) = 
\sum_{k \in {\ZZ}_+} \sum_{s \in J(k)} c_{k,s}(f)M_{k,s}, 
\end{equation}  
where we recall that $J(k) := \ \{s \in \ZZ:\ 0 \le s  \le   2^k \}$ 
is the set of $s$ for which $M_{k,s}$ do not vanish identically on  $\II$.
From the equality $M_{k,s}(2^{-k}s')= \delta_{s,s'}$ one can see that $Q_k(f)$ 
interpolates $s$ at the dyadic points $2^{-k}s, \ s \in J(k)$, i.e.
\begin{equation} \label{eq:InterpolationUnivariateQ_k(f)}
Q_k(f,2^{-k}s)= \ f(2^{-k}s), \ s \in J(k).
\end{equation}
Because of the interpolation property \eqref{eq: Interpolation[Pi_k]} and 
\eqref{eq:InterpolationUnivariateQ_k(f)}, the operators $\Pi_k$ and $Q_k$ are
interpolants. Therefore, the representations \eqref{eq:Representation[r=1]} and 
\eqref{eq:B-splineRepresentation[r=2]}
are interpolant representations. We will see that the interpolant representation 
\eqref{eq:B-splineRepresentation[r=2]} coincides  
with the classical Faber-Schauder series. The univariate 
Faber-Schauder system of functions is defined by
\begin{equation*} 
\Ff := \{\varphi^2_{k,s}: s\in Z_2(k),\ k \in {\ZZ}_+\}, 
\end{equation*}
where $Z_2(0) := \{0,1\}$ and  $Z_2(k) := \{0,1,..., 2^{k-1} - 1\}$ for $k > 0$,
\begin{equation*} 
\varphi^2_{0,0}(x):= M_{0,0}(x), \ \varphi^2_{0,1}(x):= M_{0,1}(x), \ x \in \II,
\end{equation*}
(an alternative choice is $\varphi_{0,1}(x):= 1$), and  for $k > 0$ and 
$s \in Z(k)$
\begin{equation*} 
\varphi^2_{k,s}(x):=  M_{k,2s+1}(x), \ x \in \II.
\end{equation*}
It is known that $\Ff$ is a basis in $C(\II)$. (See \cite{KS} for details about the Faber-Schauder system.)

By a direct computation we have for the component functions $q_k(f)$ in the piecewise linear
quasi-interpolant representation \eqref{eq:B-splineRepresentation[r=2]}:
\begin{equation} \label{eq:[q_k(f),r=2]}
q_k(f)
\  = \
\sum_{s \in Z_2(k)}\lambda^2_{k,s}(f) \varphi^2_{k,s}(x).
\end{equation}
where
\begin{equation*} 
\lambda^2_{k,s}(f) \ := \
- \frac {1}{2} \Delta_{2^{-k}}^2 f(2^{-k + 1}s),   \, k > 0, \ 
\text{and} \ \lambda^2_{0,s}(f) \ := \ f(s). 
\end{equation*}

Hence, the interpolant representation \eqref{eq:B-splineRepresentation[r=2]} 
can be rewritten as the Faber-Schauder series:
\begin{equation*} 
f \ = \sum_{k \in {\ZZ}_+} \ q_k(f) = 
\sum_{k \in {\ZZ}_+} \sum_{s \in Z_2(k)} \lambda^2_{k,s}(f) \varphi^2_{k,s}, 
\end{equation*}
and for any continuous function $f$ on $\II$, 
\begin{equation*} 
\|f - Q_k(f)\|_\infty 
\ \le \
\omega_2(f, 2^{-k})_\infty.  
 \end{equation*}

Put $Z_r^d(k):= \prod_{i=1}^d Z_r(k_i), \ r=1,2$. 
For $k \in {\ZZ}^d_+$, $s \in Z_r^d(k)$, define 
\begin{equation*} 
\varphi^r_{k,s}(x)
\ := \
\prod_{i=1}^d \varphi^r_{k_i,s_i}(x_i),
\end{equation*}
and $\lambda^r_{k,s}(f) $ in the manner of the definition 
\eqref{def:Mixed[a_{k,s}(f)]} by
\begin{equation*}   
\lambda^r_{k,s}(f) 
\ := \   
\lambda^r_{k_1,s_1}((\lambda^r_{k_2,s_2}(... \lambda^r_{k_d,s_d}(f))).
\end{equation*}

\begin{theorem} \label{Theorem:F-SchRepresentation}
Let $r=1,2$, $\ 0 < p, \theta \le \infty$ and $1/p < \alpha < r$. Then there hold the following assertions. 
\begin{itemize}
\item[(i)] A function 
$f \in \MB$ can be represented by the series 
\begin{equation} \label{eq:F-SchRepresentation}
f \ =  
\sum_{k \in {\ZZ}^d_+} \sum_{s \in Z_r^d(k)} \lambda^r_{k,s}(f)\varphi^r_{k,s}, 
\end{equation}
converging in the quasi-norm of $\MB$.
Moreover, we have 
\begin{equation*} 
B^*(f)
\ := \
\left(\sum_{k \in {\ZZ}^d_+}  \left\{2^{(\alpha - 1/p)|k|_1}
\left(\sum_{s \in Z_r^d(k)}|\lambda^r_{k,s}(f)\}|^p \right)^{1/p} \right\}^\theta\right)^{1/\theta} 
 \ \le \ 
C B(f).
\end{equation*}  
\item[(ii)]  If in addition, $r=2$ and $\alpha < \min (2,1 + 1/p)$,
then a continuous function $f$ on ${\II}^d$ belongs to 
the Besov space $\MB$ if and only if 
$f$ can be represented by the series \eqref{eq:F-SchRepresentation}.
Moreover, the Besov quasi-norm $B(f)$ is equivalent to the discrete quasi-norm $B^*(f)$.
\end{itemize} 
\end{theorem}

\begin{proof}
If $r=2$, from the definition and \eqref{eq:[q_k(f),r=2]}
we can derive that 
for functions on $\II^d$ and 
$k \in \ZZ^d_+$, the component function $q_k(f)$ in the interpolant 
representation \eqref{eq:B-splineRepresentation} related to the interpolant \eqref{eq:NodalQ_k},
can be rewritten as
\begin{equation} \label{eq:F-SRepresentationq_k(f)}
q_k(f)
\  = \
\sum_{s \in Z_2^d
(k)}\lambda^2_{k,s}(f) \varphi^2_{k,s}(x).
\end{equation}  
Therefore, Theorem \ref{Theorem:F-SchRepresentation} 
is the rewritten Theorem \ref{Theorem:Representation}. 
This does not hold for the case $r=1$. However,
the last case can be proven in a way completely similar to the proof of Theorem \ref{Theorem:Representation}
by using the above mentioned properties of the functions $\varphi^1_{k,s}$ and operators $\Pi_k$.
\end{proof}

For $m \in {\ZZ}_+$, we have by \eqref{eq:F-SRepresentationq_k(f)}
\begin{equation*}
R^r_m(f) 
\ = \ 
R_m(f) 
\ = \ 
\sum_{k \in \Delta(m)} \sum_{s \in Z_r^d(k)} \lambda^r_{k,s}(f)\varphi^r_{k,s}.
\end{equation*} 
 For functions $f$ on ${\II}^d$, $R^r_m$ defines a linear sampling algorithm 
of the form \eqref{def:L_n(f)} on the grid $G_r^d(m)$ where 
$G_r^d(m):= \{ 2^{-k}s: k \in \Delta(m),\ s \in I_r^d(k)\}$,
$I_1^d(k):=  \{s \in {\ZZ}^d_+: 0 \le s_i \le 2^{k_i} - 1, \ i \in N[d]\}$,  $I_2^d(k) := I^d(k)$. 
More precisely, 
\begin{equation*} 
R^r_m(f) 
\ = \ 
L^r_n(f) 
\ = \ 
\sum_{k \in \Delta(m)} \sum_{j \in I_r^d(k)} f(2^{-k}j) \psi^r_{k,j}, 
\end{equation*} 
where
\begin{equation*} 
n \ := \ 
\sum_{k \in \Delta(m)} |I_r^d(k)|
\ \asymp \ 
2^m  m^{d-1};
\end{equation*}
\begin{equation*} 
\psi^r_{k,s}(x)
\ = \
\prod_{i=1}^d \psi^r_{k_i,s_i}(x_i),\ k \in {\ZZ}^d_+, \ s \in I_r^d(k), 
\end{equation*}
and the univariate functions $\psi^r_{k,s}$ are defined by
\begin{equation*} 
\psi^1_{k,s}
\ =  \ 
\begin{cases}
\varphi^1_{k,s}, & \ k=0, \ s=0, \\
\varphi^1_{k,j}, & \ k>0, \ s = 2j+1, \\
- \varphi^1_{k,j}, & \ k>0, \ s=2j,
\end{cases}
\end{equation*}
and
\begin{equation*} 
\psi^2_{k,s}
\ =  \ 
\begin{cases}
\varphi^2_{k,s}, & \ k=0, \\
- \frac{1}{2} \varphi^2_{k,0}, & \ k>0, \ s=0, \\
\varphi^2_{k,j}, & \ k>0, \ s = 2j+1, \\
- \frac{1}{2} ( \varphi^2_{k,j} + \varphi^2_{k,j-1}), & \ k>0, \ s=2j, \\
- \frac{1}{2} \varphi^2_{k,2^{k-1}-1}, & \ k>0, \ s= 2^k.
\end{cases}
\end{equation*}
 
From the interpolation properties \eqref{eq: Interpolation[Pi_k]} and \eqref{eq:InterpolationUnivariateQ_k(f)},  
the equality $\varphi^r_{k,s}(2^{-k}s')= \delta_{s,s'}$ 
one can easily verify that $R^r_m(f)$ interpolates $f$ at the grid $G_r^d(m)$, i.e., 
\begin{equation*} 
R^r_m(f,x)= \ f(x), \ x \in G_r^d(m).
\end{equation*}

\begin{theorem} \label{Theorem:LowerBoundE(m),r=1}
Let $r= 2$, $\ 0 < p, q, \theta \le \infty$,  and $1/p < \alpha < \min (2,1 + 1/p)$. Then we have 
\begin{itemize} 
\item[{\rm (i)}] For $p \ge q$,
\begin{equation*} 
E(m)
 \ \gg \ 
2^{- \alpha m} m^{(d-1)(1 - 1/\theta)_+}.
\end{equation*}
\item[{\rm (ii)}] For $p < q$, 
\begin{equation*} 
E(m)
 \ \gg \ 
2^{- (\alpha - 1/p + 1/q) m} m^{(d-1)(1/q - 1/\theta)_+}.
\end{equation*}
\end{itemize} 
\end{theorem}

\begin{proof} 
Put $\Gamma (m):= \{k \in {\NN}^d: |k|_1 = m + 1 \}$. 
Let the half-open $d$-cube $I(k,s)$ be defined by
$I(k,s):= \prod_{i=1^d} [s_i 2^{-(k_i - 1)}, (s_i + 1)2^{-(k_i - 1)})$.  
Notice that $I(k,s) \subset {\II}^d$ and $I(k,s) \cap I_(k,s') = \emptyset$ for $s \ne s'$. 
Moreover,  if $0 < \nu \le \infty$, for $k \in \Gamma (m), \ s \in Z^d(k)$,
\begin{equation} \label{asymp:Normvarphi_{k,s}}
\|\varphi^2_{k,s}\|_\nu 
\ = \ 
\left(\int_{I(k,s)}|\varphi^2_{k,s}(x)|^\nu  dx \right)^{1/\nu}
\ \asymp \ 2^{- m/\nu} , 
\end{equation}
with the change to sup when $\nu = \infty$, and 
\begin{equation}  \label{asymp:NormSumvarphi_{k,s}}
\left \|\sum_{s \in Z_2^d(k)}\varphi^2_{k,s}\right\|_\nu 
\ \asymp \ 1. 
\end{equation}

\noindent
{\it Case} (i):
For an integer $m \ge 1$, we take the functions
\begin{equation} \label{def:g_1}
g_1 
\ := \ 
C_1 2^{-\alpha m} \sum_{s \in Z_2^d({\bar k})} \varphi^2_{{\bar k}, s} 
\end{equation}
with some ${\bar k} \in \Gamma (m)$, and
\begin{equation} \label{def:g_2}
g_2 
\ := \ 
C_2 2^{- \alpha m} m^{-(d-1)/\theta} \sum_{k \in \Gamma (m)} \sum_{s \in Z^d(k)} \varphi^2_{k,s}.
\end{equation}
Notice that the right side of \eqref{def:g_1} and \eqref{def:g_2}
 defines the series \eqref{eq:F-SchRepresentation} of $g_i, \ i=1,2$.
By Theorem \ref{Theorem:F-SchRepresentation}
and \eqref{asymp:NormSumvarphi_{k,s}}
we can choose constants $C_i$ so that 
$g_i \in B^\alpha_{p,\theta}$ for all $m \ge 1$ and $i=1,2$. It 
is easy to verify that $g_i - R^2_m(g_i) = g_i \ i=1,2$. 
We have by \eqref{asymp:NormSumvarphi_{k,s}} 
\begin{equation*}
E(m)
\  \ge \ \|g_1\|_q  
\  \gg \ 
 2^{-\alpha m} 
\end{equation*}
if $\theta \le 1$, and 
\begin{equation*}
E(m)
\  \ge \ \|g_2\|_q \ \ge \|g_2\|_{q^*}
\ \gg \ 
  2^{-\alpha m} m^{(d-1)(1 - 1/\theta)}
\end{equation*}
if $\theta > 1$, where $q^*:= \min(q,1)$.

\noindent
{\it Case} (ii): 
Let $s(k)\in {\ZZ}^d_+$ be defined by $s(k)_i = \sum_{j=1}^{k_i-2} 2^j$ if $k_i > 2$, 
and $s(k)_i = 0$ if $k_i = 2$ for $i = 1,...,d$,
 and $\Gamma^* (m):= \{k \in \Gamma (m):  k_i \ge 2, \ i = 1,...,d\}$. 
For an integer $m \ge 2$, we take the functions
\begin{equation} \label{def:g_3}
g_3 
\ = \ 
C_3 2^{-(\alpha - 1/p)m} \varphi^2_{k^*,s(k^*)} 
\end{equation}
with some $k^* \in \Gamma^* (m)$, and
\begin{equation} \label{def:g_4}
g_4 
\ = \ 
C_4 2^{-(\alpha - 1/p)m} m^{-(d-1)/\theta} \sum_{k \in \Gamma^* (m)} \varphi^2_{k,s(k)}.
\end{equation}
Similarly to the functions $g_i, \ i =1,2$, the right side of \eqref{def:g_3} and \eqref{def:g_4}
 defines the series \eqref{eq:F-SchRepresentation} of $g_i, \ i =3,4$, and we can choose constants $C_i$ so that 
$g_i \in B^\alpha_{p,\theta}$ for all $m \ge 2$ and $i=3,4$. Obviously,  $g_i - R^2_m(g_i) = g_i, \ i=3,4$. 
We have by \eqref{asymp:Normvarphi_{k,s}}
\begin{equation*}
E(m)
\  \ge \ \|g_3\|_q  
\  \gg \ 
 2^{-(\alpha - 1/p + 1/q)m} 
\end{equation*}
if $\theta \ge q$, and
\begin{equation*}
E(m)
\  \ge \ \|g_4\|_q 
\  \gg \ 
 2^{-(\alpha - 1/p + 1/q)m} m^{(d-1)(1/q - 1/\theta)} 
\end{equation*}
if $\theta < q$.
\end{proof}

From Theorems \ref{Theorem:UpperBoundE(m),alpha>1/p} and \ref{Theorem:LowerBoundE(m),r=1} 
we obtain

\begin{theorem} \label{Theorem:AsympE(m),r=1} 
Let  $r= 2$, $\ 0 < p, q, \theta \le \infty$, and $1/p < \alpha < \min (2,  1 + 1/p)$. Then we have 
\begin{itemize} 
\item[{\rm (i)}] For $p \ge q$,
\begin{equation*} 
E(m)
 \ \asymp \ 
\begin{cases}
2^{- \alpha m}, \ & \theta \le \min(q,1), \\
2^{- \alpha m} m^{(d-1)(1 - 1/\theta)}, \ & \theta > 1, q \ge 1.
\end{cases}
\end{equation*}
\item[{\rm (ii)}] For $p < q <\infty$, 
\begin{equation*} 
E(m)
 \ \asymp \ 
2^{- (\alpha - 1/p + 1/q) m} m^{(d-1)(1/q - 1/\theta)_+}.
\end{equation*}
\end{itemize} 
\end{theorem}

Notice that Theorem \ref{Theorem:AsympE(m),r=1}(i)
has been proven in \cite{SU} for the 
$1 \le p= q =\theta \le \infty$.  

\section{Appendix}
\label{Appendix}

\begin{lemma} \label{Lemma: Asymp[lambda_n]}
Let $1 \le p, q \le \infty$, $\ 0 <  \theta \le \infty$ 
and $\alpha > (1/p - 1/q)_+$. Then we have the following asymptotic order of 
$\lambda_n(B^\alpha_{p,\theta})_q$. 
\begin{itemize}
\item[{\rm (i)}] For $p \ge q$,
\begin{equation*} 
\lambda_n(B^\alpha_{p,\theta})_q
 \ \asymp \
\begin{cases}
(n^{-1} \log^{d-1}n)^\alpha, \ &  \theta \le 2 \le q \le p < \infty, \\
(n^{-1} \log^{d-1}n)^\alpha, \ & \theta \le 1, \ p = q = \infty, \\
(n^{-1} \log^{d-1}n)^\alpha, \ & 1 < p = q \le 2, \  \theta \le q,  \\
(n^{-1} \log^{d-1}n)^\alpha(\log^{d-1}n)^{1/q - 1/\theta}, \ & 1 < p = q \le 2, \ \theta > q \\
(n^{-1} \log^{d-1}n)^\alpha (\log^{d-1}n)^{ 1/2 - 1/\theta}, \ & \theta > 2.
\end{cases}
\end{equation*}
\item[{\rm (ii)}] For $1 < p < q < \infty$, 
\begin{equation*} 
\lambda_n(B^\alpha_{p,\theta})_q
 \ \asymp \ 
 \begin{cases}
 (n^{-1} \log^{d-1}n)^{\alpha - 1/p + 1/q}, \ & 2 \le p, \ 2 \le \theta \le q,  \\
(n^{-1} \log^{d-1}n)^{\alpha - 1/p + 1/q}(\log^{d-1}n)^{(1/q - 1/\theta)_+}, \ &  q \le 2.
\end{cases}
 \end{equation*}
\end{itemize} 
\end{lemma}

\begin{proof}
This lemma was proven in \cite{G}, \cite{R1} except the cases 
$\theta \le 2 \le q \le p < \infty$ and $\theta \le 1, \ p = q = \infty$ which can be obtained from
the asymptotic order \cite{R1}
\begin{equation*} 
\lambda_n(B^\alpha_{p,\theta})_q
 \ \asymp \
(n^{-1} \log^{d-1}n)^\alpha,
\begin{cases}
 & 1 \le \theta \le 2 \le q \le p < \infty, \\
 & \theta = 1, \ p = q = \infty,
\end{cases}
\end{equation*}
and the equalities $\lambda_n(W)_q = \lambda_n(\text{co} W)_q$ and 
$\text{co} B^\alpha_{p,\theta} = B^\alpha_{p,\max(\theta,1)}$, where $\text{co} W$ denotes the convex 
hull of $W$. 
\end{proof}

For ${\bf p} = (p_1,...,p_d) \in (0,\infty)^d$, we defined the mixed integral  quasi-norm 
$\|\cdot\|_{\bf p}$ for functions on ${\II}^d$ as follows
\begin{equation*}
\|f\|_{\bf p}
\ := \
\left( \int_{\II} \left(\cdots \int_{\II}\left(\int_{\II}|f(x)|^{p_1} dx_1\right)^{p_2/p_1} 
dx_2 \cdots \right)^{p_d/p_{d-1}} dx_d \right)^{1/p_d},
\end{equation*}
and put $1/{\bf p} := (1/p_1,...,1/p_d)$.
If ${\bf p},{\bf q} \in (0,\infty)^d$ and ${\bf p} \le {\bf q}$, then there holds 
Nikol'skii's inequality for any $f \in \Sigma_r^d(k)$,
\begin{equation} \label{ineq:Nikol'skii}
\|f\|_{\bf q}
\  \le \
C 2^{|(1/{\bf p} - 1/{\bf q})k|_1} \|f\|_{\bf p}
\end{equation}
with constant $C$ depending on ${\bf p},{\bf q},d$ only. This inequality can be proven by  a generalization of the
Jensen's inequality for mixed norms and the following equivalences of the 
mixed integral  quasi-norm $\|\cdot\|_{\bf p}$. 
For all $k \in {\ZZ}^d_+$  and all $f \in \Sigma_r^d(k)$ of the form \eqref{def:StabIneq}, 
\begin{equation*} 
 \|f\|_{\bf p} 
\ \asymp \ \prod_{i=1}^d 2^{- k_i/p_i} \|\{a_s\}\|_{{\bf p},k}, 
\end{equation*}
where
\begin{equation*}
\| \{a_s\} \|_{{\bf p},k}
 \ := \
\left( \sum_{s_d \in J(k_d)} \left(\cdots \sum_{s_2 \in J(k_2)}\left(\sum_{s_1 \in J(k_1)}
|a_s|^{p_1}\right)^{p_2/p_1}  \cdots \right)^{p_d/p_{d-1}} \right)^{1/p_d}.
\end{equation*}

\begin{lemma}  \label{Lemma:IneqInt[varphi_k*varphi_s]}
Let $0 < p < q < \infty$, $\delta = 1/2 - p/(p+q)$. 
If $k,s \in {\ZZ}^d_+$, then for any $\varphi_k \in \Sigma_r^d(k)$ and 
$\varphi_s \in \Sigma_r^d(s)$, there holds the inequality 
\begin{equation*}  
\int_{{\II}^d} |\varphi_k(x) \varphi_s(x)|^{q/2} dx
\ \le \ 
C A_k A_s 2^{- \delta |k - s|_1},  
\end{equation*}
with some constant $C$ depending at most on $p,q,d$, where
\begin{equation*}  
A_k := \ \left(2^{(1/p - 1/q)|k|_1} \|\varphi_k \|_p\right)^{q/2}.
\end{equation*} 
\end{lemma}

\begin{proof} 
 Put $\nu := (p + q)/p$. Then $\delta = 1/2 - 1/\nu$ and $2 < \nu < \infty$.
Let $\nu'$ be given by $1/\nu + 1/\nu' = 1$. Then $1 < \nu' < 2 $.
Let ${\bf u}, {\bf v} \in (0,\infty)^d$ be defined by ${\bf u} := q{\bf v}/2$ and 
$v_i = \nu$ if $k_i \ge s_i$ and $v_i = \nu'$ if $k_i < s_i$ for $i= 1,...,d$.  
 Let ${\bf u}'$ and ${\bf v}'$ be given by $1/{\bf u}+ 1/{\bf u}' = {\bf 1}$
and $1/{\bf v}+ 1/{\bf v}' = {\bf 1}$, respectively.
Notice that ${\bf v} \in (1,\infty)^d$.
Applying H\"older's inequality for the mixed norm $\|\cdot\|_{\bf v}$
to the functions $|\varphi_k|^{q/2}$ and $|\varphi_s|^{q/2}$, we obtain 
\begin{equation}  \label{ineq:Int[varphi_k*varphi_s]}
\int_{{\II}^d} |\varphi_k(x) \varphi_s(x)|^{q/2} dx
\ \le \ 
\||\varphi_k|^{q/2}\|_{\bf v} \||\varphi_s|^{q/2}\|_{{\bf v}'}
\ = \ 
\|\varphi_k\|_{\bf u}^{q/2} \|\varphi_s\|_{{\bf u}'}^{q/2}.
\end{equation}
Since ${\bf u} > p{\bf 1}$ and ${\bf u}' > p{\bf 1}$, by the inequality \eqref{ineq:Nikol'skii} we have
\begin{equation} \label{ineq:Norms[varphi_k,varphi_s]}
\|\varphi_k\|_{\bf u}
\  \le \
2^{|({\bf 1}/p - 1/{\bf u})k|_1} \|\varphi_k\|_p, \quad
 \|\varphi_s\|_{{\bf u}'}
 \  \le \
2^{|({\bf 1}/p - 1/{\bf u}')s|_1} \|\varphi_s\|_p.  
\end{equation}
From \eqref{ineq:Int[varphi_k*varphi_s]} and 
\eqref{ineq:Norms[varphi_k,varphi_s]} we prove the lemma.
\end{proof}

\begin{lemma}  \label{Lemma:IneqL_qNorm<L_pNorm}
Let $0 < p < q < \infty$ and  
$g \in L_q$ be represented by  the  series 
\begin{equation*} 
g \ = \sum_{k \in {\ZZ}^d_+} \ g_k, \ g_k \in \Sigma_r^d(k).
\end{equation*}
Then there holds the inequality 
\begin{equation} \label{ineq:IneqL_qNorm<L_pNorm}
\|g\|_q
 \ \le C \left( \sum_{k \in {\ZZ}^d_+} \ \|2^{(1/p - 1/q)|k|_1} g_k \|_p^q \right)^{1/q}, 
\end{equation}
with some constant $C$ depending at most on $p,d$, whenever the right side is finite. 
\end{lemma}

\begin{proof} 
It is enough to prove the inequality \eqref{ineq:IneqL_qNorm<L_pNorm} for $g$ of the form
\begin{equation*} 
g \ = \sum_{k \le m} \ g_k, \ g_k \in \Sigma_r^d(k),
\end{equation*} 
for any $m \in {\ZZ}^d_+$. 

Put $n:= [q] + 1$. Then $0 < q/n \le 1$. By Jensen's inequality we have
\begin{equation*}
\begin{aligned}  
\left|\sum_{k \le m} \ g_k(x)\right|^q
\ & = \
\left(\left|\sum_{k \le m} \ g_k(x)\right|^{q/n}\right)^n \\
  \  \le \
\left(\sum_{k \le m} \ |g_k(x)|^{q/n}\right)^n 
 \ & = \
\sum_{k^1 \le m} \cdots  \sum_{k^n \le m}\ \prod_{j=1}^n |g_{k^j}(x)|^{q/n}. \\  
\end{aligned}
\end{equation*}
consequently,
\begin{equation} \label{ineq:[|f|_q^q](1)}
\|g\|_q^q
  \  \le \
\sum_{k^1 \le m} \cdots  \sum_{k^n \le m}\ \int_{{\II}^d} \ \prod_{j=1}^n |g_{k^j}(x)|^{q/n} dx.   
\end{equation}
By use of the identity
\begin{equation} \label{ineq:prod_{j=1}^n a_j}
\prod_{j=1}^n a_j 
\ = \
\left( \prod_{i \ne j} a_i a_j \right)^{1/2(n-1)} 
\end{equation}
for non-negative numbers $a_1,..., a_n$, we get
\begin{equation*} 
J:= \ \int_{{\II}^d} \ \prod_{j=1}^n |g_{k^j}(x)|^{q/n} dx 
\ = \ 
\int_{{\II}^d} \prod_{i \ne j}|g_{k^i}(x)g_{k^j}(x)|^{q/2n(n-1)} dx.   
\end{equation*}
Hence, applying H\"older's inequality to $n(n-1)$ functions 
in the right side of the last equality, Lemma \ref{Lemma:IneqInt[varphi_k*varphi_s]} and 
\eqref{ineq:prod_{j=1}^n a_j} gives
\begin{equation*}
\begin{aligned}  
J
\ & \le \ 
\prod_{i \ne j} \left( \int_{{\II}^d} |g_{k^i}(x)g_{k^j}(x)|^{q/2} dx \right)^{1/n(n-1)} 
\  \le \ 
\prod_{i \ne j} \left( A_{k^i} A_{k^j} 2^{- \delta |k^i - k^j|_1} \right)^{1/n(n-1)} \\   
\ & = \ 
\prod_{i \ne j} \left( A_{k^i} A_{k^j} \right)^{1/n(n-1)}
\left\{\left( \prod_{i \ne j}\prod_{i' = 1}^n 2^{- \delta |k^i - k^{i'}|_1}
\prod_{j' = 1}^n 2^{- \delta |k^j - k^{j'}|_1}\right)^{1/2(n-1)}\right\}^{1/n(n-1)} \\  
\ & = \ 
\left\{\prod_{i \ne j} \left( A_{k^i} A_{k^j} \right)^{2/n}
\left( \prod_{i' = 1}^n 2^{- \delta |k^i - k^{i'}|_1}
\prod_{j' = 1}^n 2^{- \delta |k^j - k^{j'}|_1}\right)^{1/n(n-1)}\right\}^{1/2(n-1)} \\
\ & = \ 
\prod_{j = 1}^n A_{k^j}^{2/n}
\left(\prod_{i = 1}^n 2^{- \delta |k^j - k^i |_1}\right)^{1/n(n-1)} 
\  = \ 
\left(\prod_{j = 1}^n A_{k^j}^2
\prod_{i = 1}^n 2^{- \lambda \delta |k^j - k^i |_1}\right)^{1/n},  
\end{aligned}  
\end{equation*}
where $\lambda := \delta/(n-1) > 0$.
Therefore, from \eqref{ineq:[|f|_q^q](1)} and H\"older's inequality we obtain
\begin{equation} \label{ineq:[|f|_q^q](2)} 
\begin{aligned}  
\|g\|_q^q
  \ & \le \
\sum_{k^1 \le m} \cdots  \sum_{k^n \le m}\ 
\left(\prod_{j = 1}^n A_{k^j}^2
\prod_{i = 1}^n 2^{- \lambda \delta |k^j - k^i |_1}\right)^{1/n} \\
  \ & \le \
\prod_{j = 1}^n 
\left(\sum_{k^1 \le m} \cdots  \sum_{k^n \le m}\ 
A_{k^j}^2
\prod_{j = 1}^n 2^{- \lambda \delta |k^j - k^i |_1}\right)^{1/n} 
\ =: \ \prod_{j = 1}^n B_j.   
\end{aligned}  
\end{equation}
We have
\begin{equation*}
\begin{aligned}  
B_j 
 \ & = \
 \sum_{k^j \le m} \ A_{k^j}^2 
\sum_{k^1 \le m} \cdots  \sum_{k^{j-1} \le m} \ \sum_{k^{j+1} \le m} \cdots  \sum_{k^n \le m} 
\prod_{i = 1}^n 2^{- \lambda \delta |k^j - k^i |_1} \\
 \ & = \
 \sum_{k^j \le m} \ A_{k^j}^2 \left( \sum_{s \le m}  
 2^{- \lambda \delta |k^j - s|_1} \right)^{n-1}
 \ \le \
 C \sum_{k^j \le m} \ A_{k^j}^2.
\end{aligned}  
\end{equation*}
Using this estimate for $B_j$, we can continue \eqref{ineq:[|f|_q^q](2)} and finish the estimation 
of $\|g\|_q^q$ as follows.
\begin{equation*}
\begin{aligned}  
\|g\|_q^q 
 \ & \le \
 \prod_{j = 1}^n B_j^{1/n} 
 \ \le \
 C \sum_{k \le m} \ A_k^2 \\
\ & = \
C \sum_{k \le m} \ \|2^{(1/p - 1/q)|k|_1} g_k \|_p^q .
\end{aligned}  
\end{equation*}
Thus, the proof of the lemma is completed.
\end{proof}

\noindent
{\bf Remark} A trigonometric polynomial version of Lemma \ref{Lemma:IneqL_qNorm<L_pNorm}
was proven in \cite{Te1b} for $1 \le p < q <\infty$.

\noindent
{\bf Acknowledgments.} This work is supported by Vietnam National Foundation for Science and 
Technology Development (NAFOSTED).

\end{document}